\RequirePackage{fix-cm}
\documentclass[smallcondensed,envcountsect]{svjour3}       % onecolumn (second format)
\smartqed  % flush right qed marks, e.g. at end of proof
 \journalname{}
\usepackage{color}
\usepackage[english]{babel}
\usepackage{verbatim}
\usepackage{xcolor}
\usepackage{graphicx}
\usepackage{multirow}
\usepackage{mathtools,amssymb}
\usepackage{numberbysection}
\usepackage{booktabs}
\usepackage{comment}
\usepackage{nicefrac}
\usepackage{longtable}
\usepackage{amssymb}
\usepackage{float}

\numberbysection

\begin{document}

\title{A simulation-based optimization approach for the calibration of a discrete event simulation model of an emergency department}

\titlerunning{Simulation-based optimization for ED DES model calibration}        % if too long for running head

%\author{Alberto De Santis~\href{https://orcid.org/0000-0001-5175-4951}{\includegraphics[scale=0.5]{orcid.png}} 
%        \and
%        Tommaso Giovannelli~\href{https://orcid.org/0000-0002-1436-5348}{\includegraphics[scale=0.5]{orcid.png}}
%        \and
%        Stefano Lucidi
%        \and
%        Mauro Messedaglia
%        \and
%        Massimo Roma~\href{https://orcid.org/0000-0002-9858-3616}{\includegraphics[scale=0.5]{orcid.png}}
%}

\author{Alberto~De~Santis~\href{https://orcid.org/0000-0001-5175-4951}{\includegraphics[scale=0.5]{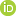}}
	\and 
	Tommaso~Giovannelli~\href{https://orcid.org/0000-0002-1436-5348}{\includegraphics[scale=0.5]{orcid.png}}
	\and 
	Stefano~Lucidi~\href{https://orcid.org/0000-0003-4356-7958}{\includegraphics[scale=0.5]{orcid.png}}
	\and
	Mauro~Messedaglia
	\and
	Massimo~Roma~\href{https://orcid.org/0000-0002-9858-3616}{\includegraphics[scale=0.5]{orcid.png}}
}

\authorrunning{A. De Santis, T. Giovannelli, S. Lucidi, M. Messedaglia, M. Roma} % if too long for running head

\institute{Alberto De Santis, Stefano Lucidi, Tommaso Giovannelli, Massimo Roma \at
                Dipartimento di Ingegneria Informatica Automatica e Gestionale ``A. Ruberti'', SAPIENZA Universit\`a di Roma, via Ariosto, 25 -- 00185 Roma, Italy. \\
              \email{desantis@diag.uniroma1.it, lucidi@diag.uniroma1.it, giovannelli@diag.uniroma1.it, roma@diag.uniroma1.it}           %  \\
           \and
           Mauro Messedaglia \at
              ACTOR Start up of SAPIENZA Universit\`a di Roma,
              via Nizza 45, 00198 Roma, Italy.
              \email{mauro.messedaglia@gmail.com}
}

\date{}

\maketitle

\begin{abstract}
Accurate modeling of the patient flow within an Emergency Department (ED) is required by all studies dealing with the increasing and well-known problem of overcrowding. Since Discrete Event Simulation (DES) models are often adopted with the aim of assessing solutions for reducing the impact of this worldwide phenomenon, an accurate estimation of the service time of the ED processes is necessary to guarantee the reliability of the results. However, simulation models concerning EDs are frequently affected by data quality problems, thus requiring a proper estimation of the missing parameters.

In this paper, a simulation-based optimization approach is used to estimate the incomplete data in the patient flow within an ED by adopting a model calibration procedure. The objective function of the resulting minimization problem represents the deviation between simulation output and real data, while the constraints ensure that the response of the simulation is sufficiently accurate according to the precision required. Data from a real case study related to a big ED in Italy is used to test the effectiveness of the proposed approach. The experimental results show that the model calibration allows recovering the missing parameters, thus leading to an accurate DES model.

%Include keywords, PACS and mathematical
%subject classification numbers as needed.
\keywords{Simulation-based Optimization \and Emergency Department \and Discrete Event Simulation \and Model Calibration \and Derivative-Free Optimization}
\end{abstract}

\section{Introduction}\label{sec:intro}
Over the last few years, Emergency Departments (EDs) have been raising an increasing attention in the Operations Research and Management Science communities due to the international phenomenon of the overcrowding \cite{hoot.aronsky:08,Paul.2010,Wiler:2011,disomma.2014,Nahhas.2017}, which leads to longer waiting times, higher mortality rates, and lower patient satisfaction \cite{pines.2008}. One of the most popular tools adopted to study this crucial problem is Discrete Event Simulation (DES) \cite{Zeinali.2015,Joshi.2016,Wong.2016,Bedoya-Valencia:2016,kramer.2020,fava.2021}, which is used to represent the complex and stochastic patient flows in the ED in place of analytical models (for a recent literature review of simulation modeling in EDs, see \cite{Salmon:2018}). Frequently, DES models are combined with an optimization algorithm to make optimal decisions in relation to some Key Performance Indicators (KPIs). The resulting Simulation-Based Optimization (SBO) approach \cite{fu.2015,gosavi.2014,amaran.2014}, whose effective application relies on the accuracy of the simulation model, has been frequently applied to hospital EDs (for a recent review, see, e.g., \cite{yousefi.2020}).     

Many papers in the literature deal with the quality of the input data used in ED simulation models. This issue, which strongly affects the reliability of the results, has been carefully analyzed. Indeed, the cost, time, and challenges required by collecting ED empirical data represents a serious limitation for every simulation model \cite{Paul.2010}. In order to accurately replicate and predict the patient flows within the ED, \cite{duma.2018} develops a process mining approach which handles the noise factors in the dataset after introducing assumptions on how to interpret the data. In particular, these factors include the following cases: starting and ending time of each activity performed in the ED may be unknown; information about urgent patients may be registered after the activity is completed and, in general, the timestamps of each activity may not be promptly recorded at the right time; for each patient, the same activity may be recorded multiple times for technical reasons, giving rise to misleading information. A framework to categorize all the ED data quality issues is proposed in \cite{vanbrabant:2019}, which also provides assessment techniques for each data quality problem category. Moreover, this paper highlights that most of the works in the literature focus on the problem of missing data, which is the problem addressed in this paper. It consists in the lack of information on some or all the key timestamps that define the activities performed in the ED, i.e., starting and ending time. This well-known issue prevents gaining knowledge about the duration of each activity, which is required to estimate the corresponding probability distributions to use in the simulation model. Several SBO approaches are proposed in the literature to tackle this crucial problem. The idea behind each approach is to leverage the known information in order to estimate the parameters of the probability distributions underlying the missing data. This is accomplished by comparing the Key Performance Indicators (KPIs) computed through the simulation model with the corresponding values derived by the data collected in the real system. The resulting procedure is known as \emph{model calibration}.

The mathematical formulation of the optimization problem and the specific algorithm used for determining the optimal solution are the two features that distinguish the papers dealing with missing data. For instance, both \cite{Kuo.2016} and \cite{Guo.2016} adopt similar approaches that leverage the time differences between the known timestamps. The former proposes an unconstrained optimization problem where the objective function is a consistency measure that compares the average, the standard deviation, and the proportions of the time differences for each triage tag. The latter uses a constrained optimization problem where the objective function is based on a modified chi--square goodness of fit and the constraints are introduced to guarantee each time difference the same level of accuracy. As regards the approaches used to solve the problem, both papers use metaheuristic procedures. In particular, \cite{Kuo.2016} considers both a descent method and a simulated annealing algorithm, while \cite{Guo.2016} adopts an approach that combines a genetic algorithm with simulated annealing and optimal computing budget allocation. Moreover, in \cite{Guo.2016} the authors point out that the approach in \cite{Kuo.2016} could be improved in several ways: including in the objective function all the time differences defined by the available timestamps, since delays in one activity may impact on downstream activities in the patient flow; modeling the possibility for each patient to have more than one medical visit, which may affect the time differences; using a more formal objective function and solving the resulting optimization problem through a more efficient algorithm. 

Another paper in the literature on missing data that is worth being mentioned is \cite{Liu.2017}, which assumes an agent--based simulation model, as opposed to the DES model considered in \cite{Kuo.2016} and \cite{Guo.2016}. Although the framework underlying this work is different, a simulation-based optimization problem is used for the same goal of minimizing the deviation between real data and simulation output in order to estimate the missing parameters of the simulation model. The Length Of Stay (LOS), i.e., the difference between the discharge time and the arrival time to the ED, is the time difference considered in the objective function, which adopts the Jensen--Shannon divergence. A systematic method based on the pattern search method APPSPACK \cite{gray.2006} is used to find the optimal configuration of parameters. 

The approach proposed in this paper aims to handle the problem of missing timestamps, which affects many simulation models dealing with EDs, as evidenced by the papers in the specific literature. Among the sources of noise affecting the quality of the input data necessary to build a reliable simulation model, the problem of missing timestamps has a strong impact on the overall accuracy of the simulation, since it prevents the knowledge of the values used to derive appropriate probability distributions. Therefore, while other issues may be resolved by either carefully cleaning the dataset or introducing assumptions on how to interpret the data, the unavailability of timestamps requires more sophisticated procedures.

Compared to the other papers dealing with missing data in ED datasets, the approach discussed in this paper aims both to propose a new formulation of the resulting optimization problem and to improve on the optimization strategies typically used in the literature. Since building a simulation model is a process that requires a considerable amount of effort and thus is not expected to be completed in a short time, an exact algorithm providing optimal solutions may be preferable to metaheuristic procedures, whose final solutions are returned faster but without optimality guarantees. However, although global convergent algorithms appear to be a reasonable choice, metaheuristics are the methods mainly adopted in the literature dedicated to missing data (see, e.g., \cite{Kuo.2016,Guo.2016}). To fill this gap, in this paper a SBO approach is developed both to propose an alternative version of the optimization problems used in \cite{Kuo.2016,Guo.2016} and to use a Derivative-Free Optimization (DFO) strategy based on global convergence, which allows the algorithm to find an optimal solution with optimality guarantees. Similarly to \cite{Kuo.2016} and differently from \cite{Guo.2016}, Weibull distributions are adopted to generate the values of the activity service times associated with missing data. Indeed, this probability distribution is considered suitable when data is unavailable (see, e.g., \cite{law:15}). Moreover, the most critical patients, who are excluded from the two approaches developed by \cite{Kuo.2016} and \cite{Guo.2016}, are included to prevent the simulation model from returning inaccurate results for this important class of patients. Finally, in the proposed constrained optimization problem, the objective function is based on the comparison between real and simulated Empirical Cumulative Distribution Functions (ECDFs), which do not rely on intervals, unlike the functions adopted in \cite{Kuo.2016,Guo.2016,Liu.2017}. 

This paper is organized as follows. Section~\ref{data_collection_ED} reports the typical structure of the datasets related to the patient flow within an ED, highlighting the available information and the KPIs of interest. Section~\ref{sec:SBO_calibration} defines the SBO problem used to calibrate a general DES model concerning an ED. Section~\ref{chap:EDPoliclinicoUmbertoI} describes the real ED used as a case study and the DES model implemented. Section~\ref{sec:model_calibration_policlinico} reports the results obtained by applying the SBO approach to the case study. Finally, in Section~\ref{sec:conc} we draw some concluding remarks and we discuss directions for future work.

\section{Data collection in ED}\label{data_collection_ED}
\begin{figure}[htbp]
	\centering
	\includegraphics[width=12truecm]{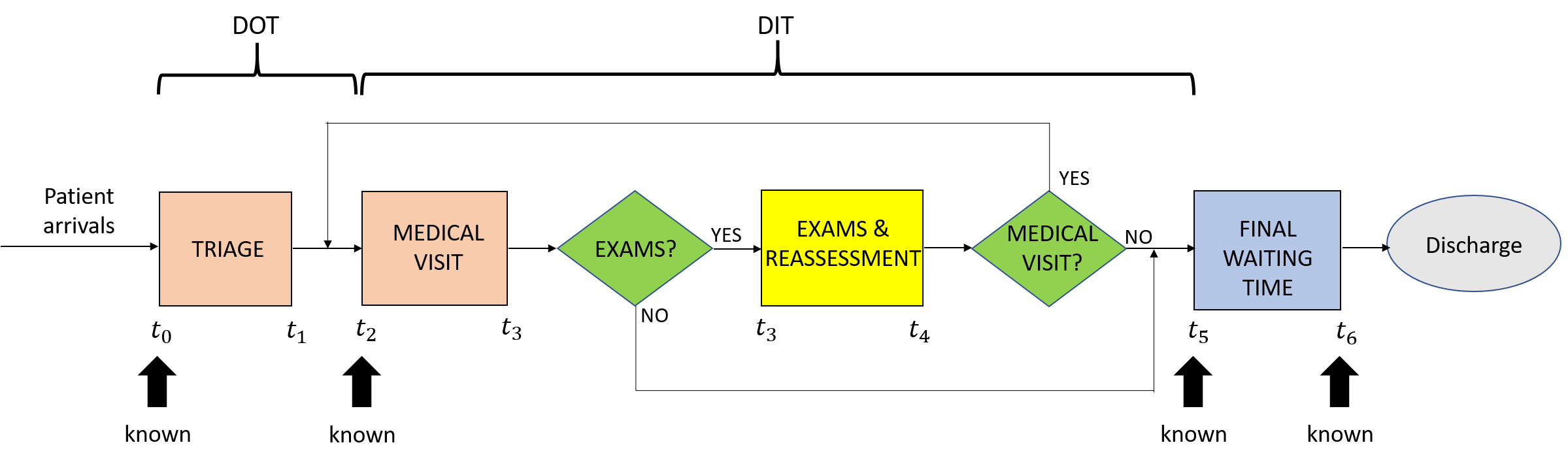}
	\caption{Timestamps collected throughout the patient flow.}\label{fig:missing_data}
\end{figure}

\noindent
The timestamps defining the activities in the patient flow are represented in Figure \ref{fig:missing_data} and described in Table \ref{tab:timestamp}.
\begin{table}[htbp]
	\small\sf\centering
	\caption{Description of the timestamps collected throughout the patient flow.} \label{tab:timestamp}
	\begin{tabular}{ll}
		\toprule
		{\sc \emph{$t_0$}} &  starting time of triage. \\
		%\midrule
		{\sc \emph{$t_1$}} & ending time of triage. \\
		%\midrule
		{\sc \emph{$t_2$}} & starting time of medical visit. \\
		%\midrule
		{\sc \emph{$t_3$}} & ending time of medical visit and starting time of examinations. \\
		{\sc \emph{$t_4$}} & ending time of the last exam or reassessment. \\
		{\sc \emph{$t_5$}} & patient receives the last medical report. \\
		{\sc \emph{$t_6$}} & patient is discharged and leaves the ED. \\
		\bottomrule
	\end{tabular}
\end{table}
%	\begin{description}
%		\item[]\emph{$t_0$}: starting time of triage;
%		\item[]\emph{$t_1$}: ending time of triage;
%		\item[]\emph{$t_2$}: starting time of medical visit;
%		\item[]\emph{$t_3$}: ending time of medical visit and starting time of examinations;
%		\item[]\emph{$t_4$}: patient receives the examination response;
%		\item[]\emph{$t_5$}: patient receives the last examination response;
%		\item[]\emph{$t_6$}: patient is discharged and leaves the hospital.		
%	\end{description}
Although these timestamps are commonly recorded by the electronic systems used to collect data, the actual patient flow may include extra times omitted from this scheme. For example, since the arrival time to the ED is usually not registered, before triage there may be an additional waiting time not included in any record. Moreover, extra waiting times and observation periods may be present both in the examination part and in the discharge phase. Since usually there is no available data for these times, in Figure \ref{fig:missing_data} both waiting and service times for exams are included in a single box between $t_3$ and $t_4$, which are associated with the end of the medical visit and the end of the last exam or reassessment after a visit, respectively. After receiving the medical report, patients may be either subject to additional medical visits or discharged after a final waiting time, which is due to a further observation or to the time required by the ED staff to prepare the discharge. It is important to point out that the disposition, which refers to the decision to admit a patient to a hospital ward, is associated with $t_6$, like the discharge. This implies that the boarding time of the admitted patients (i.e., the waiting time before being transferred to the hospital ward) is assumed to start at $t_6$, thus being excluded from the ED scope. This choice is in accordance with the case study considered in Section \ref{chap:EDPoliclinicoUmbertoI}. Contrarily, in some EDs the discharge of patients admitted may be considered as the time of admission to the hospital ward (see, e.g., \cite{Baumlin.2010}). 

In most cases, the available timestamps are $t_0$, $t_2$, and $t_6$, which are marked in the scheme as {\em known}. Since this is a typical setting in practice, throughout this paper they are supposed to be the only timestamps recorded by the ED along with $t_5$ (which may not be always available). As a consequence, for each patient it is possible to compute
\begin{itemize}
	\item the \textit{Door-to-Doctor Time} (\textbf{DOT}), which is the time difference between the start of the triage and the start of the medical visit, namely $t_2-t_0$;
	\item the \textit{Doctor-to-Discharge Time} (\textbf{DIT}), which is the time difference between the start of the medical visit and the discharge, namely  $t_5-t_2$ if $t_5$ is available, $t_6-t_2$ otherwise.
\end{itemize}
It is important to remark that when patients require additional medical visits after the exams, $DIT$ can be interpreted as the time difference between the timestamp of the last medical report (or the discharge if $t_5$ is not available) and the starting time of the first medical visit. Note that $t_4$ and $t_5$ are equal if a patient does not need additional medical visits. 

Since in many cases only a subset of the timestamps related to the patient flow are known \cite{vanbrabant:2019}, the service time cannot be computed for all the activities. Apart from the final waiting time, which can be recovered for each patient through the difference $t_6 - t_5$, note that the durations of all the other activities are not completely defined. In fact, in case of triage and medical visit, the ending time is missing, while for exams both timestamps are unknown.        

\section{Simulation-based optimization problem for ED model calibration}
\label{sec:SBO_calibration}
\noindent
The goal of the approach proposed is to recover the information needed to build an accurate simulation model by leveraging the known information through the minimization of the deviation between real data and simulation output.
In order to define the simulation-based optimization problem, let us introduce the sets below.
\begin{itemize}
	\vspace{-0.2cm}
	\item Let $C$ be the set of the triage tags.
	\item Let $U(c)$ be the set of units where patients with tag $c \in C$ can be visited and treated.
	\item Let $\mathcal{T}=\{DOT,DIT\}$ be the set of the time differences considered.
\end{itemize}
In the absence of data, suitable probability distributions are Weibull and lognormal (see \cite{law:15}). Let us arbitrarily consider a Weibull distribution, whose probability density function is $$ f(w) = \begin{cases}
	\alpha\beta^{-\alpha}w^{\alpha-1}e^{-(w/\beta)^\alpha} \text{ if } w > 0, \\ 0 \text{ otherwise,} 
\end{cases}$$ where $\alpha > 0$ and $\beta > 0$ are the shape and the scale parameters, respectively. Therefore, a Weibull distribution is assumed for the triage, medical visit, and exams service times and, for each of them, different pairs of shape and scale parameters are considered based on the triage tag $c$ and on the unit $u$. This choice leads to a number of different pairs of parameters equal to $\sum_{c \in C}|U(c)|$, where $|U(c)|$ refers to the cardinality of the set $U(c)$. Let us denote as $x \in \mathbb{R}^{n_1}$, $y \in \mathbb{R}^{n_2}$, and $z \in \mathbb{R}^{n_3}$ the corresponding probability distribution parameters for the service times considered, where $n_v \le 2 \, \sum_{c \in C}|U(c)|$ with $v \in \{1,2,3\}$. In particular, for all $c \in C$ and for all $u \in U(c)$, the shape and scale parameters of the Weibull distributions are 
\begin{itemize}
	\item $x_1^{c \, u}$ and $x_2^{c \, u}$ for the triage probability distribution,
	\item $y_1^{c \, u}$ and $y_2^{c \, u}$ for the medical visit probability distribution,
	\item $z_1^{c \, u}$ and $z_2^{c \, u}$ for the exams probability distribution.
\end{itemize}  
Let $F^{sim}_{c \, u \, i}$ and $F^{real}_{c \, u \, i}$ be the ECDFs of the values of the simulated and real time difference $i \in \mathcal{T}$ for patients with tag $c$ visited in unit $u$. Moreover, let $k^{sim}_{c \, u \, i}$ and $k^{real}_{c \, u \, i}$ be the number of such patients from the simulation and from the real dataset. Hence, for all $j \in \{sim,real\}$, we have that
$$F^j_{c \, u \, i}(t) = \frac{1}{k^j_{c \, u \, i}}\sum_{h=1}^{k^j_{c \, u \, i}} \textbf{1}_{\{\tau_{{c \, u \, i} \,h} \le t\}} \quad \text{ with } t \ge 0,$$
where $\tau_{{c \, u \, i} \, h}$ is the value of the time difference $i$ recorded for the $h$--patient, with $h \in \{1,\ldots,k^j_{c \, u \, i}\}$, considered from $j \in \{sim,real\}$. It is important to point out that the values $\tau_{{c \, u \, i} \,h}$ of the time differences computed from the simulation depend on the service times of triage, medical visit, and exams drawn from the Weibull distributions described above. Therefore, this dependence can be written explicitly as $F^{sim}_{c \, u \, i}(t;x,y,z)$, where $x$, $y$, and $z$ are the vectors containing all the associated shape and scale parameters.  

The mathematical problem formulation is reported as follows
\vspace{-0.5cm}
\begin{center}
	\begin{equation} \label{prob:model_calibration}
		\begin{split}
			\min_{x, \, y, \, z}\text{ } &\sum_{c \in C}\sum_{u \in U(c)}\sum_{i\in\mathcal{T}} \Big(\int_{0}^{\infty} (F^{sim}_{c \, u \, i}(t;x,y,z)-F^{real}_{c \, u \, i}(t))^2 \, dt \Big) \\
			s.t. \text{ } \, & x\in {\cal P},
		\end{split}
	\end{equation} 	
\end{center}
where the feasible set
%	\begin{equation*} 
%	\begin{array}{l}
%	{\cal P}=\Bigl\{x\in{\mathbb{R}}^{n_1}, y\in{\mathbb{R}}^{n_2}, z\in{\mathbb{R}}^{n_3} ~ | ~  g_i(x,y,z)\leq 0 \text{ and } h_i(x,y,z)\leq 0 \bigr. \cr
%	\ \cr
%	\bigl.\hspace{1truecm} \text{ for all } c \in C, \, u \in U(c), \, i \in \mathcal{T}_{c \, u} \Bigr\}
%	\end{array}
%	\end{equation*}
\begin{equation*} 
	\begin{array}{l}
		{\cal P}=\Bigl\{(x,y,z) \in {\mathbb{R}}^{n_1} \times {\mathbb{R}}^{n_2} \times {\mathbb{R}}^{n_3} ~ | ~ \vspace{0.2truecm} \\ \vspace{0.2truecm} \qquad \qquad \qquad g_{c \, u \, i}(x,y,z)\leq 0, \\ \vspace{0.2truecm} \qquad \qquad \qquad h_{c \, u \, i}(x,y,z)\leq 0,  \\ \vspace{0.2truecm} \qquad \qquad \qquad l_x \le x \le u_x \\ \vspace{0.2truecm} \qquad \qquad \qquad l_y \le y \le u_y \\ \vspace{0.2truecm} \qquad \qquad \qquad l_z \le z \le u_z \bigr. \cr
		\ \cr
		\bigl.\hspace{0.7truecm} \text{ for all } c \in C, \, u \in U(c), \, i \in \mathcal{T}  \Bigr\}
	\end{array}
\end{equation*}	
is defined by the following functions
\begin{itemize}
	\item $g_{c \, u \, i}(x,y,z) = \big|\frac{\mu^{sim}_{c \, u \, i}(x,y,z) - \mu^{real}_{c \, u \, i}}{\mu^{real}_{c \, u \, i}}\big| - tol_{\mu}^{c \, u \, i}$, which compares the sample means $\mu^{sim}_{c \, u \, i}(x,y,z)$ and $\mu^{real}_{c \, u \, i}$ of the time difference $i$ computed through the simulated (by averaging over the independent replications) and real data, respectively,
	\item $h_{c \, u \, i}(x,y,z) = \big|\frac{\sigma^{sim}_{c \, u \, i}(x,y,z) - \sigma^{real}_{c \, u \, i}}{\sigma^{real}_{c \, u \, i}}\big| - tol_{\sigma}^{c \, u \, i}$, which compares the sample standard deviations $\sigma^{sim}_{c \, u \, i}(x,y,z)$ and $\sigma^{real}_{c \, u \, i}$ of the time difference $i$ computed through the simulated (by averaging over the independent replications) and real data, respectively,
\end{itemize}
and $l_x, l_y, l_z, u_x, u_y$, and $u_z$ are vectors defining bound constraints. Positive thresholds, namely $tol_{\mu}^{c \, u \, i}$ and $tol_{\sigma}^{c \, u \, i}$, are used to state the degree of accuracy required for the simulation model.
The objective function is the sum of the integrals of the squared difference between $F^{sim}_{c \, u \, i}$ and $F^{real}_{c \, u \, i}$ over the sets $C$, $U(c)$, and $\mathcal{T}$. The choice of using ECDFs instead of histograms, which are commonly adopted in the literature \cite{Kuo.2016,Guo.2016}, is due to the different reliability of the information provided. In particular, compared to ECDFs, the description of the data provided by histograms is strongly affected by the choice of the width of each bin, which is not straightforward and may lead to distributions with different shapes based on how data is grouped. The decision variables are the Weibull distribution parameters contained in the vectors $x, y$, and $z$. For the sake of simplicity, the previous formulation does not include variables that are not parameters of probability distributions, although they may be present. It is important to remark that in the general framework described, the dependence of each pair of parameters on both $c$ and $u$ implies that the triage service time is affected by the triage tag and the ED unit. This assumption, which turns out to be reasonable when applied to the duration of medical visit and exams, may lead to an excessive number of variables for the triage service time, if it does not hold. Indeed, while the color tag significantly affects the triage duration since urgent patients undergo a faster triage than less critical patients, the impact of the unit may be negligible. However, as concerns the case study described in Section~\ref{chap:EDPoliclinicoUmbertoI}, interviews with the ED staff have shown that slightly different procedures are adopted by the nurses in charge of triage based on the unit where a patient is assigned. This means that final conclusions should be drawn after analyzing the specific system considered. Since the probability distribution parameters affect the dimension of the optimization problem, avoiding unnecessary variables allows the algorithm to benefit from a lower computational cost.	

\section{Case study: ED of Policlinico Umberto I in Rome}
\label{chap:EDPoliclinicoUmbertoI}

\noindent
The case study concerns the ED of {\em Policlinico Umberto~I}, which is a very large hospital in Rome, Italy. It is the biggest ED in the region of Lazio in terms of number of patient arrivals per year (about 140,000 on average). By using the data collected from the patient flow through the ED for the whole year 2018, this case study is adopted to test the effectiveness of the approach proposed in this paper. 
%gli arrivi nei 4 mesi per ciascun giorno della settimana sono, partendo dal lunedì:3431+3084+3149+3081+3253+3042+2830
\par
\subsection{Description of the ED} \label{sec:description_ED_Policlinico}
\noindent
The ED of Policlinico Umberto I is divided into several areas, each one associated with a medical specialty. The backbone of the ED is the \emph{central area}, which is devoted to treating diseases and disorders related to internal medicine and general surgery, which affect the majority of patients. Separated from this main area, there are other parts of the ED that deal with the following medical specialties: ophthalmology, obstetrics, pediatrics, hematology, and dentistry. The focus of this paper is on the central area, whose number of arrivals in 2018 amounted to more than 50,000.  
\par
In order to gain a complete understanding of the ED processes, the description of the ED units and staff and the summary of the patient flow are reported in the sequel. Other than a \emph{triage area}, where each incoming patient is assigned a color tag by a nurse in charge of this task, the ED is composed by
\begin{itemize}
	\item[$\bullet$] a \emph{Medical Unit (MU)}, devoted to patients needing specialized medications and treatments, with areas dedicated to the most critical patients;
	\item[$\bullet$] a \emph{Surgical Unit (SU)}, devoted to patients needing either to receive a surgical operation or to recover from it, with areas dedicated to the most critical patients;
	\item[$\bullet$] a \emph{Resuscitation Area (RA)}, for the most acutely ill and injured patients, who need timely treatments;		
	\item[$\bullet$] a \emph{Minor Injuries Unit (MIU)}, for the least urgent patients, whose treatment can be delayed or deferred;
	\item[$\bullet$] an \emph{Orthopedic Unit (OU)}, for patients suffering from orthopedic disorders.
\end{itemize}
\noindent
Moreover, all of these units have rooms where patients can either wait for exams or stay for observation. Red-tagged patients can be visited and treated in RA or in dedicated areas within MU and SU, which are open 24 hours a day and are provided with equipment and staff specialized for dealing with life-threatening illnesses and injuries. In particular, 1 and 2 seats are available in the dedicated areas of MU and SU, respectively, and further 2 seats are available in RA. As concerns the medical treatment of the other patients, MU and SU can host up to three and two patients during the day (8.00 a.m.--8.00 p.m.), respectively. At night, MU can host two patients, while in SU one seat is available. Moreover, MIU has two seats, which are available from 8.00 a.m. to 8.00 p.m., Monday through Saturday. When patients experience excessive waiting times, two additional seats may be used to visit up to four patients simultaneously. All this information is summarized in Tables~\ref{tab:unit_shifts}--\ref{tab:unit_shifts_red}.  
\begin{table}[htbp]
	\caption{Number of seats available for medical visit and treatment in MU, SU, and MIU, the latter being open from Monday to Saturday.}\label{tab:unit_shifts}
	\small\sf\centering
	\begin{tabular}{lccc}
		%\hline
		% after \\: \hline or
		%   \cline{2-3}
		& MU & SU & MIU  \\ \midrule
		\multicolumn{1}{l}{Day (8.00 a.m.--8.00 p.m.)} & 3 & 2 & 2 \\ %\hline			
		\multicolumn{1}{l}{Night (8.00 p.m.--8.00 a.m.)} & 2 & 1 & 0 \\
		%\midrule
	\end{tabular}
\end{table}
\begin{table}[htbp]
	\caption{Number of seats available  for medical visit and treatment of red-tagged patients in RA and in the dedicated areas of MU and SU.}\label{tab:unit_shifts_red}
	\small\sf\centering
	\begin{tabular}{lccc}
		%\hline
		% after \\: \hline or
		%   \cline{2-3}
		& MU & SU & RA  \\ \midrule
		\multicolumn{1}{l}{Day (8.00 a.m.--8.00 p.m.)} & 1 & 2 & 2 \\ %\hline			
		\multicolumn{1}{l}{Night (8.00 p.m.--8.00 a.m.)} & 1 & 2 & 2 \\
		%\midrule
	\end{tabular}
\end{table}
\begin{table}[htbp]
	\small\sf\centering
	\caption{Feasible assignments of patients to the ED units according to the color tag. A cross at the entry $(i,j)$ indicates that a patient with color tag $i$ can be assigned to the unit $j$. }\label{tab:triagetag_assignment_policlinico}
	\begin{tabular}{lccccc}
		%\toprule
		%\cline{2-5}
		%& \multicolumn{2}{|c||}{} & \multicolumn{2}{|c|}{} \\
		& MU & SU & RA & MIU & OU \\
		%& \multicolumn{2}{|c||}{} & \multicolumn{2}{|c|}{} \\
		\midrule
		\multicolumn{1}{l}{\sc White} & - & - & -  & X  & X \\ %\hline
		\multicolumn{1}{l}{\sc Green} & X & X & - & X & X  \\ %\hline
		\multicolumn{1}{l}{\sc Yellow} & X & X & - & - & X  \\ %\hline
		\multicolumn{1}{l}{\sc Red} & X & X & X & - & - \\
		%\bottomrule
	\end{tabular}
\end{table}
\par
As regards the patient flow, after arriving autonomously or by emergency medical vehicles, all the incoming patients are admitted to the triage area, where a nurse assigns the color tag. After the triage, the patients are visited and treated in one of the units previously described, according to the color tag assigned and the severity of the illness/injury. Table~\ref{tab:triagetag_assignment_policlinico} represents a scheme showing the units where patients may be assigned based on the color tags. In case of red-tagged patients, the medical visit is timely performed in RA or in the dedicated areas of MU and SU. As concerns the other color tags, the yellow and green tagged patients share the same resources in MU and SU. However, while the former patients can be assigned only to MU and SU, the latter may be sent to MIU during its opening hours if their health conditions are deemed as not likely to worsen. If MIU is closed, all the green-tagged patients are visited and treated in MU and SU. This diversion allows the ED to reduce the occurrence of work overload in SU and MU, which may give rise to overcrowded units. Moreover, it is important to point out that the white triage tag is assigned only if MIU is open, otherwise the green tag is used. 
\par
In many cases, after the medical visit, additional exams may be required. Other than performing reassessments of the patients and requiring additional exams, physicians may also require further observation periods. Finally, at the end of the pathway, patients are discharged from the ED. This last stage includes a final waiting time whose duration depends on the type of outcome. Indeed, a longer wait is expected for patients that need to either be hospitalized at a hospital ward or transferred to another hospital, while patients discharged home can usually leave the ED in shorter time. 

\subsection{The Discrete Event Simulation model}\label{sec:desm_policlinico}
\begin{figure}
	\centering
	\includegraphics[width=\textwidth]{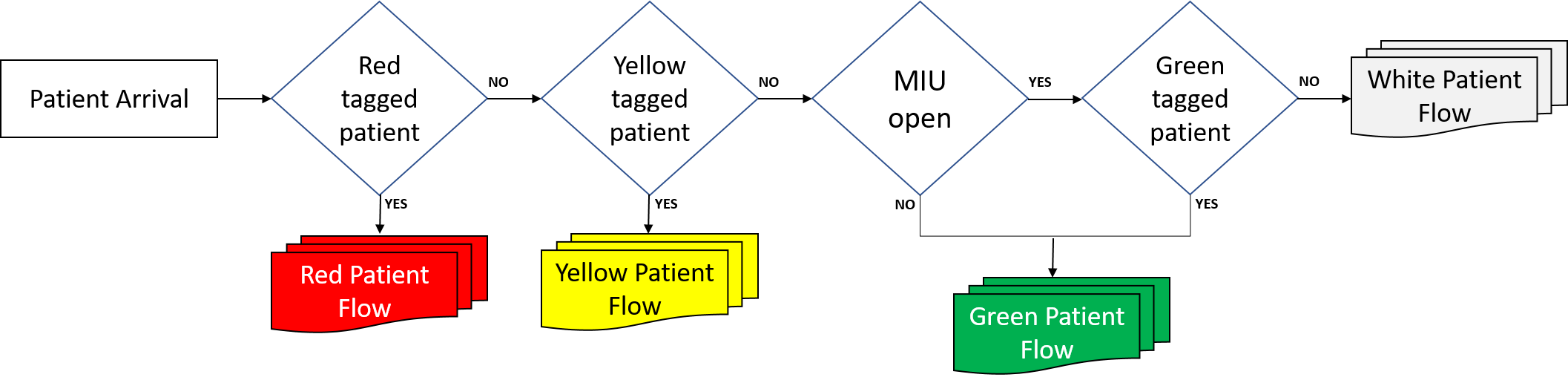}
	\caption{Plot of the patient flow  in the ED from the arrival to the assignment of the color tag.}
	\label{fig:patient_flow}
\end{figure}

\begin{figure}
	\centering
	\includegraphics[width=0.8\textwidth]{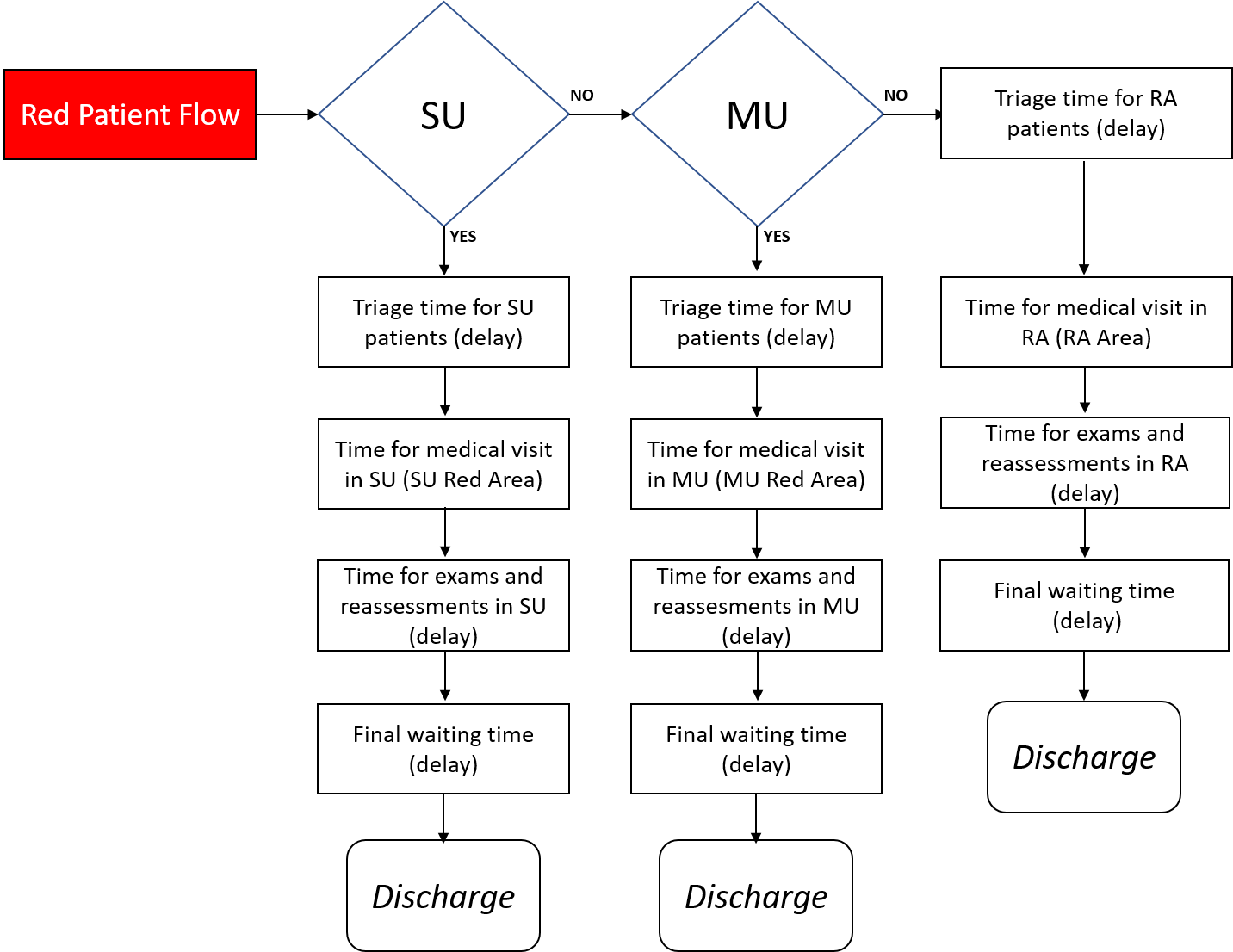}
	\caption{Plot of the red-tagged patient flow.}
	\label{fig:patient_flow_R}
\end{figure}

\begin{figure}[!tbp]
	\centering
	\begin{minipage}[b]{0.4\textwidth}
		\includegraphics[width=\textwidth]{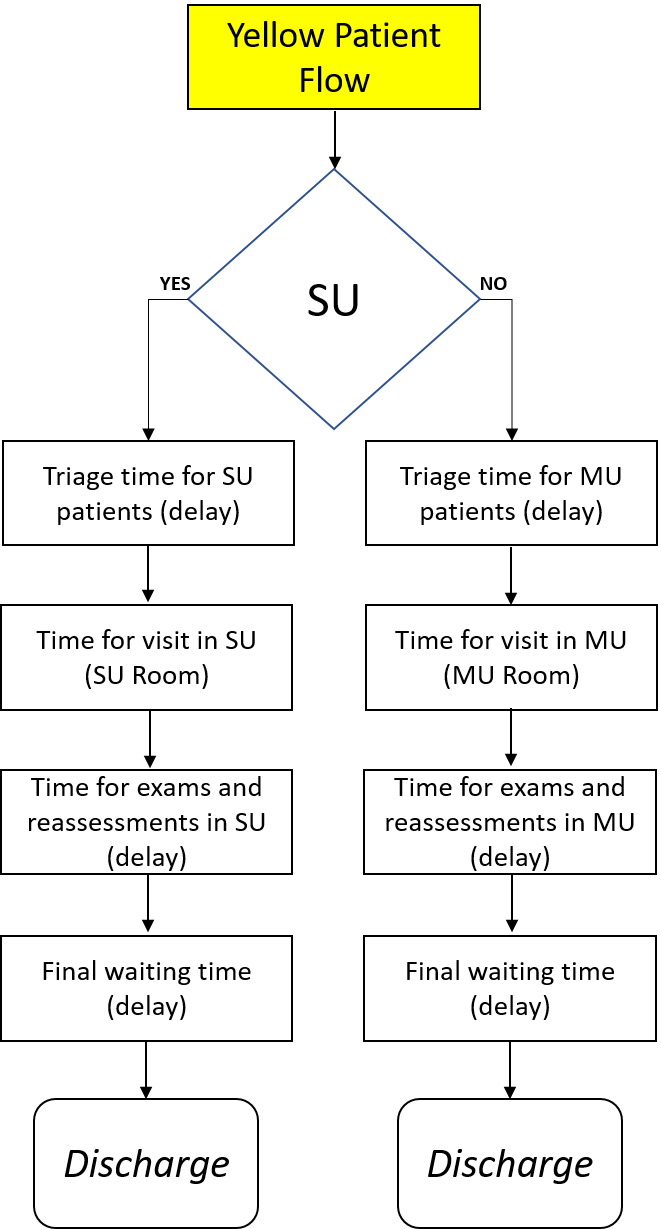}
		\caption{Plot of the yellow tagged patient flow.}\label{fig:patient_flow_Y}
	\end{minipage}
	\hfill
	\begin{minipage}[b]{0.4\textwidth}
		\includegraphics[width=\textwidth]{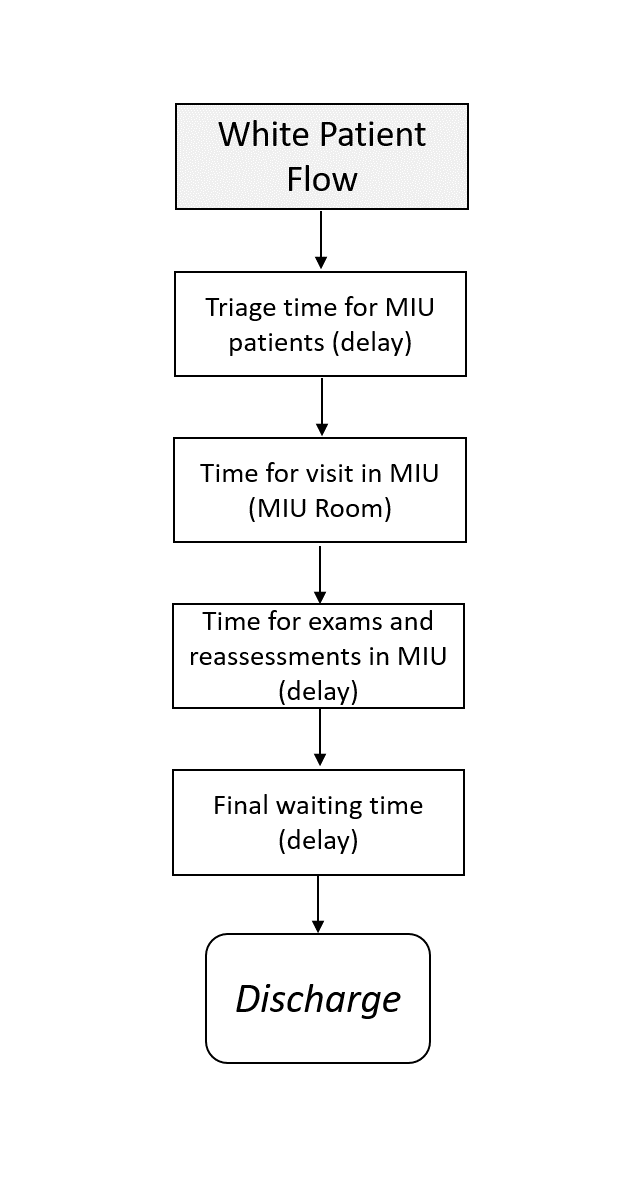}
		\caption{Plot of the white tagged patient flow.}\label{fig:patient_flow_W}
	\end{minipage}
\end{figure}

\begin{figure}
	\centering
	\includegraphics[width=0.8\textwidth]{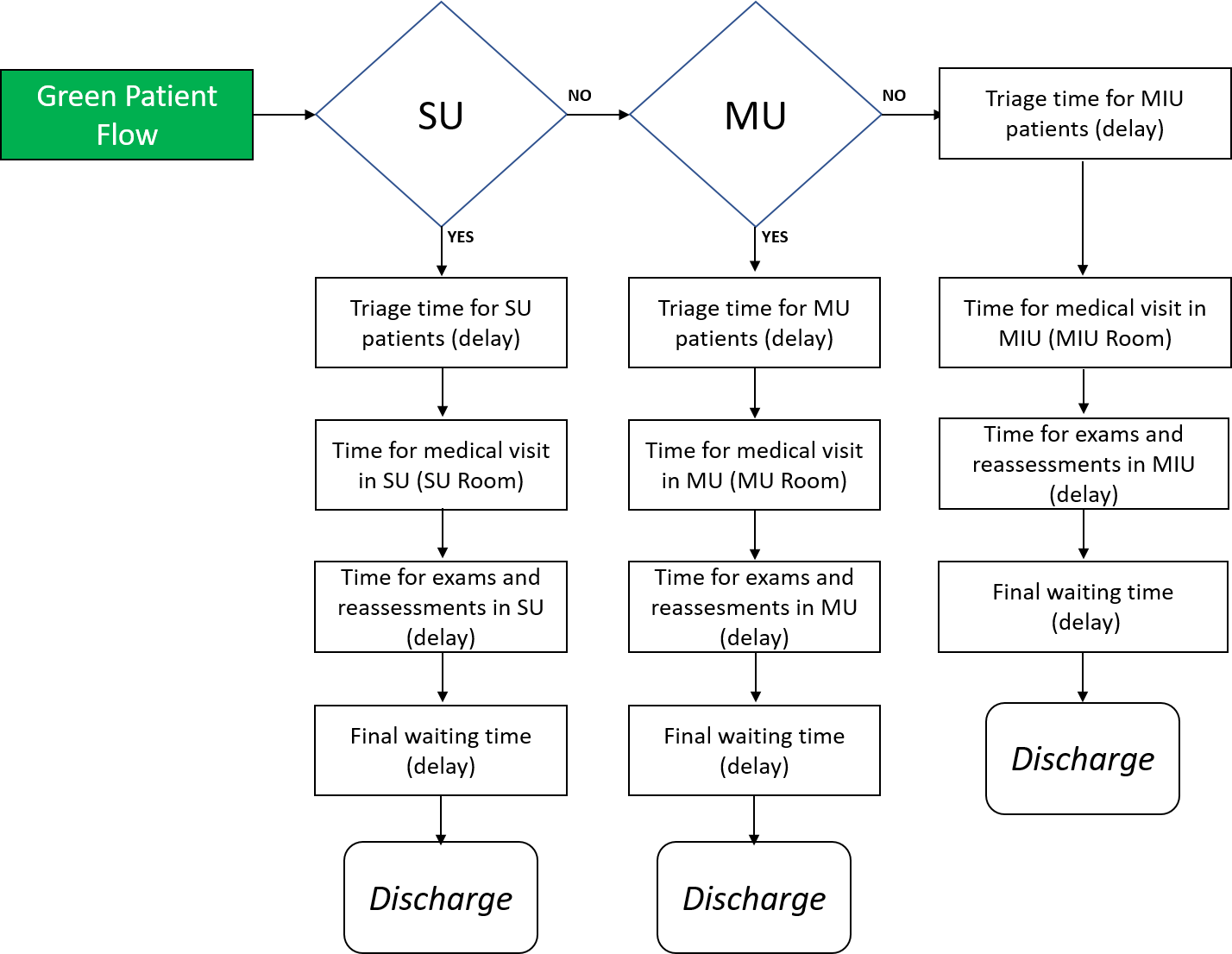}
	\caption{Plot of the green-tagged patient flow.}
	\label{fig:patient_flow_G}
\end{figure}

This section describes the DES model of the ED of Policlinico Umberto I, which has been implemented by using Simmer \cite{simmer.2019}, a process-oriented and trajectory-based DES package for {\sf R}. In the DES model, patients are represented by the model \textit{entities}, which are created according to the statistical model used for modeling the arrival process. After arriving at the ED, the simulated patient flow followed by each entity is based on trajectories determined by the logical rules used to build the model. Each trajectory is associated with the pathway followed by a patient according both to the color tag received at triage and to the ED unit assigned. Figures~\ref{fig:patient_flow}--\ref{fig:patient_flow_G} show the patient flow from the arrival to the discharge according to the color tag assigned at triage. In particular, Figure~\ref{fig:patient_flow} focuses on the logic underlying the assignment of color tag at triage. The other figures deal with the segments of the patient flow following the triage, providing in brackets the information about the \textit{resources} seized in case the corresponding process is not a simple delay. Such resources represent the seats at each ED unit, whose capacity can be either fixed, as is the case with RA and the areas of MU and SU dedicated to red patients, or based on schedule, as for MU, SU, and MIU.   
\par	
The patient flow through the model can be described as follows. After being created at the beginning of the simulation, the entities corresponding to deceased patients leave the model, while the remaining entities enter the {\em triage area} for the assignment of both color tag and ED unit, which are stored as entity attributes. The discrete uniform probability distributions used for assigning the color tag and the ED unit vary according to the time at which the entity starts the triage. If this event happens in the daytime, the corresponding probability distributions include also the white color tag and MIU among the alternatives to sample. The triage phase is represented by as many delay processes as the number of color tags and units. Each of these delays is associated with a probability distribution that returns the value of the triage duration for each patient. In particular, the 8 different processes considered correspond to the possible pairs of color tags and units reported in Table~\ref{tab:triagetag_assignment_policlinico}, with the exception of the pair given by the green color tag and MIU. Indeed, the triage duration for a green-tagged patient eligible for MIU depends on the specialty required, whether medical or surgical. Since this information is unknown, in the DES model half of the MIU patients are subject to the triage duration of MU patients, while the other half are subject to the triage duration of SU patients. It is important to remark that, since the processes are simple delay, no queue is generated before the triage. The reason underlying this choice is twofold: on the one hand, this allows for alignment between data e simulation model, since the timestamps denoted as $t_0$ in Figure~\ref{fig:missing_data} are used to derive the arrival probability distribution, which consequently returns the starting time of triage (and not the starting time of the waiting time for triage); on the other hand, the waiting time before triage can be considered negligible, according to the ED staff.  
\par 
After the triage process, some entities leave the model, thus representing the patients who leave without being seen, while the other ones start waiting for the medical visit. The entity selected for the visit depends both on the priority class, i.e., the color tag, and on the FIFO criterion. Then, when the visit begins, one seat in the corresponding unit is seized and the duration of the visit is returned by a suitable probability distribution based on the entity color tag.
\par
After the medical visit, the following phase includes {\em exams and reassessments}, whose duration is generated by means of a suitable probability distribution. Like the triage, a single delay process is used for this phase as well, due to the lack of knowledge of the resources required. It is important to point out that all possible further treatments are considered included in the service time of this process. After the exams,
% the geometric probability distributions reported in Table~\ref{tab:timesvisit_policlinico} return the number of times a patient repeats the visit and the exams. If these processes are no longer required, 
the entities corresponding to patients {\em refusing hospitalization}, {\em leaving during exams}, and {\em being transferred to another hospital} are removed from the model, while the other entities proceed to the next stage, i.e., the {\em final waiting time} before discharge.  

Since in the dataset some timestamps defining the starting and ending time of triage, medical visit, and exams are not available, the service times corresponding to these activities cannot be recovered directly. Moreover, additional delay processes are considered to reproduce all the service times that cannot be directly computed through the data, such as the setup times that are sometimes needed for sanitizing the ED areas and the idle times caused by unexpected requests for personnel from other ED units, which give rise to sudden activity disruption. In order to effectively model these times, which impact on all the patients without predictable patterns, uniform probability distributions are considered before the queue for the medical visit, according to the color tag. 
%Due to the lack of specific data, a model calibration procedure is required to determine good estimates for the parameters of the probability distributions underlying these processes, as described in Section~\ref{sec:model_calibration_policlinico}. 
\par
%\begin{figure}
%	\centering
%	\includegraphics[width=0.8\textwidth]{PatientFlowYW.png}
%	\caption{Plot of the yellow and white tagged patient flow.}
%	\label{fig:patient_flow_Y}
%\end{figure}
%\begin{figure}
%	\centering
%	\includegraphics[width=0.40\textwidth]{PatientFlowW.png}
%	\caption{Plot of the white-tagged patient flow.}
%	\label{fig:patient_flow_W}
%\end{figure} 
\par
As regards the KPIs of interest, the focus is on the time differences $DOT$ and $DIT$ described in Figure~\ref{fig:missing_data}, which can be computed through the data of the real system and then compared with the corresponding values returned by the simulation. Therefore, the KPIs of interest are
\begin{itemize}
	\item \textit{$DOT$}, which is the time difference between the starting time of the triage and the starting time of the visit, namely $t_2-t_0$ in Figure~\ref{fig:missing_data};
	\item \textit{$DIT$}, which is the time difference between the starting time of the visit and the time of the discharge, namely $t_6-t_2$ in Figure~\ref{fig:missing_data}.
\end{itemize}

\subsection{Input analysis}
\begin{figure}[htbp]
	\centering
	\includegraphics[width=12truecm,height=5truecm]{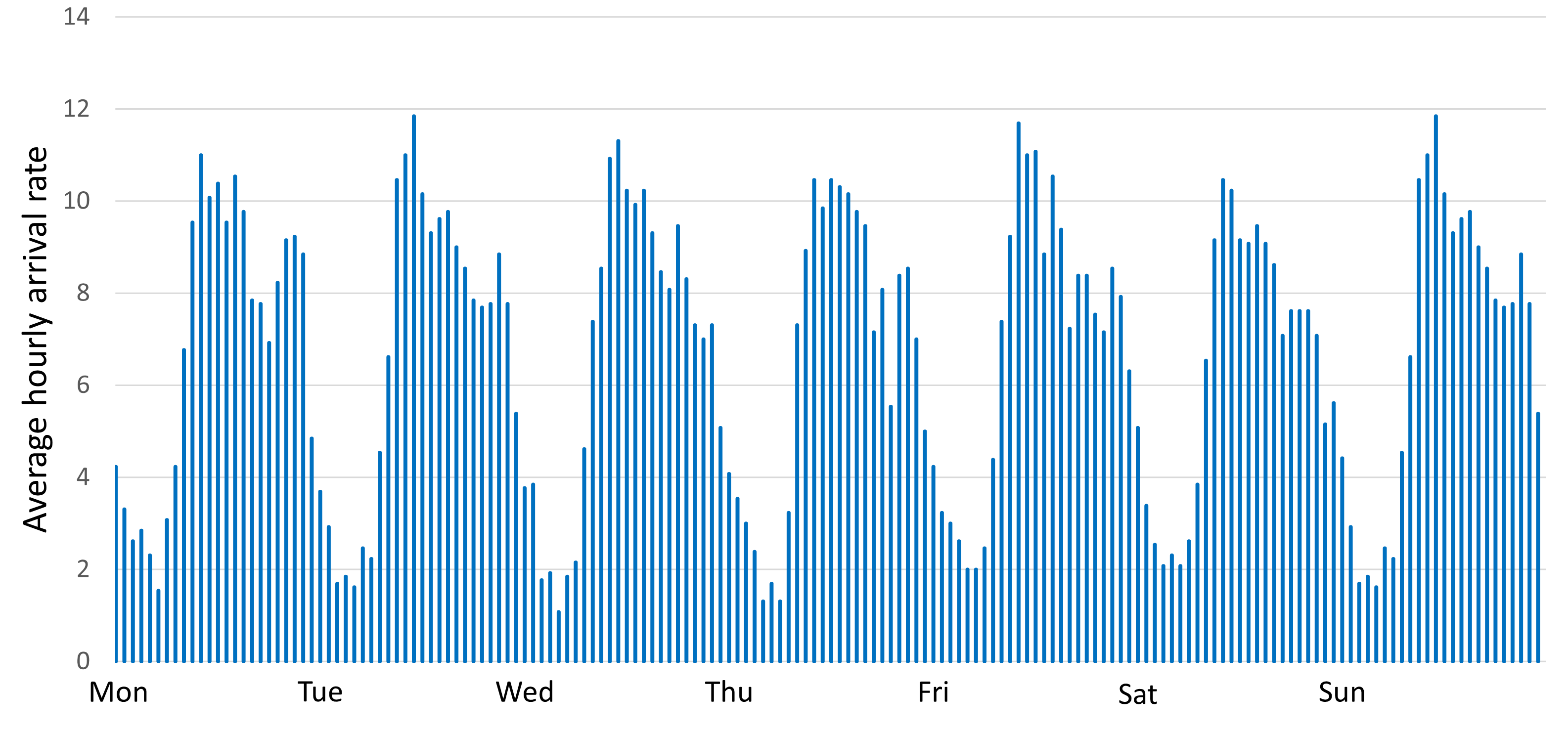}
	\caption{Plot of the weekly average hourly arrival rate for the first 13 weeks of the year.}\label{fig:weekly_arrival_rate}
\end{figure}

To analyze the arrival process, we focus on the data collected from the 1st of January to the 31st of March. In Figure \ref{fig:weekly_arrival_rate}, the weekly average hourly arrival rate obtained by averaging the number of arrivals occurring in the same hourly time slot over the 13 weeks considered is reported.
It is worth observing that, in accordance with the literature (see, i.e., \cite{Kim.2014}), the average arrival rate among the days of the week is significantly different. Therefore, since averaging over these days would lead to inaccurate results, the different days of the week must be considered separately.

The stochastic process defined by the arrivals to the ED is assumed to be well modeled by a Nonhomogeneous Poisson Process (NHPP), which is a standard assumption in the literature (see, e.g.,\cite{Whitt.2017,Kuo.2016,Zeinali.2015,Ahmed.2009,Ahalt.2018,Guo.2017}). The validity of this hypothesis for the case study considered is supported by the extensive analysis performed in \cite{DGLMR:tr,desantis.2021}, where statistical hypothesis testing is adopted. To achieve an accurate representation of the arrival rate, 24 time slots of unitary length are considered for each day of the week, thus obtaining a piecewise constant approximation. While this approach allows taking into account the within-day variation, the day-to-day variation is considered by estimating the hourly arrival rate function separately for each day of the week.

Although the collected data concerns the whole year 2018, the DES model of the ED has been built by using the data related to January. The choice of focusing on this month stems from the will of the ED management to reduce the overcrowding level observed in the winter season, which exhibits longer waiting times, as emerged by interviewing the ED staff. Among the winter months, January is observed to suffer from the heaviest workload, which puts a strain on the ED processes, thus requiring a careful analysis. However, the simulation model could be also easily adapted to include input parameters estimated from data related to different months.  
\par
From 00:00 of January 1 to 23:59 of January 31, 2018, the total number of patients arrived to the ED is 4192. The timestamps recorded are the ones marked as known in Figure~\ref{fig:missing_data}. In Table~\ref{tab:triagetag_daynight_policlinico}, the number and the percentage of color tags assigned at triage are reported along with the number of patients who leave without being seen (LWBS). Since from Monday to Saturday MIU is open from 8.00 a.m. to 8.00 p.m., the daytime and the night are analyzed separately to reflect the change in the ED setting. Although, in some cases, re-evaluations can lead to different patient tags at discharge with respect to those assigned at triage, in the dataset considered this information is unavailable. 
%\begin{table}[htbp]
%	\small\sf\centering
%	\caption{Number and percentage of color tags assigned at triage (columns 2-3) and number of patients LWBS (column 4).}\label{tab:triagetag_policlinico}
%	\begin{tabular}{lcc|c}
%		%\toprule
%		%\cline{2-5}
%		%& \multicolumn{2}{|c||}{} & \multicolumn{2}{|c|}{} \\
%		& \multicolumn{2}{c}{TRIAGE TAGS} & \multicolumn{1}{c}{LWBS}\\
%		%& \multicolumn{2}{|c||}{} & \multicolumn{2}{|c|}{} \\
%		\midrule
%		\multicolumn{1}{l}{\sc White} & 65 & 1.55\% & 1    \\ %\hline
%		\multicolumn{1}{l}{\sc Green} & 1804 & 43.17\%&  16   \\ %\hline LWBS 1.54\%
%		\multicolumn{1}{l}{\sc Yellow} & 2058 & 49.25\% & -  \\ %\hline LWBS 0.89\%
%		\multicolumn{1}{l}{\sc Red} & 252 & 6.03\% &  - \\ \midrule
%		\multicolumn{1}{l}{}  & 4179   &  100\%      &   \\
%		%\bottomrule
%	\end{tabular}
%\end{table}	
\begin{table}[htbp]
	\sf\centering
	\caption{Number and percentage of color tags assigned at triage in the daytime (columns 2-3) and at night (column 4-5).}\label{tab:triagetag_daynight_policlinico}
	\begin{tabular}{lcc|cc|c}
		%\toprule
		%\cline{2-5}
		%& \multicolumn{2}{|c||}{} & \multicolumn{2}{|c|}{} \\
		& \multicolumn{4}{c}{TAGS ASSIGNED AT TRIAGE} & \\
		& \multicolumn{2}{c}{Day (8.00 a.m.--8.00 p.m.)} & \multicolumn{2}{c}{Night (8.00 p.m.--8.00 a.m.)} & \multicolumn{1}{c}{LWBS}\\
		%& \multicolumn{2}{|c||}{} & \multicolumn{2}{|c|}{} \\
		\midrule
		\multicolumn{1}{l}{\sc White} & 65 & 2.17\% & -  & -   & 1  \\ %\hline
		\multicolumn{1}{l}{\sc Green} & 1306 & 44.32\%& 498 & 40.42 \% & 16   \\ %\hline
		\multicolumn{1}{l}{\sc Yellow} & 1420 & 48.18\% & 638 & 51.79 \% & -  \\ %\hline
		\multicolumn{1}{l}{\sc Red} & 157 & 5.33\% & 95 & 7.71 \% & - \\ \midrule
		\multicolumn{1}{l}{}  & 2948   & 100\%       & 1231  &  100\% & 17 \\
		%\bottomrule
	\end{tabular}
\end{table}
Finally, Tables~\ref{tab:triagetag_units_count_policlinico}--\ref{tab:triagetag_units_count_policlinico_V} show the number and proportion of the color tags among the ED units. Although OU is out of the scope of this analysis, Table~\ref{tab:triagetag_units_count_policlinico} includes this unit since OU patients share the triage station with the other patients, thus affecting the counts reported in Table \ref{tab:triagetag_daynight_policlinico}.   
\begin{table}[htbp]
	\scriptsize\sf\centering
	\caption{Number (and percentage) of patients assigned to the ED units for each color tag.}\label{tab:triagetag_units_count_policlinico}
	\begin{tabular}{lccccc}
		%\toprule
		%\cline{2-5}
		%& \multicolumn{2}{|c||}{} & \multicolumn{2}{|c|}{} \\
		& MU & SU & RA & MIU & OU \\
		%& \multicolumn{2}{|c||}{} & \multicolumn{2}{|c|}{} \\
		\midrule
		\multicolumn{1}{l}{\sc White} & - & - & -  & 47 (73.44 \%)  & 17 (26.56 \%) \\ %\hline
		\multicolumn{1}{l}{\sc Green} & 248 (13.75 \%) & 628 (34.81 \%) & - & 260 (14.41 \%) & 668 (37.03 \%)  \\ %\hline
		\multicolumn{1}{l}{\sc Yellow} & 1316 (63.95 \%) & 693 (33.67 \%) & - & - & 49 (2.38 \%)  \\ %\hline
		\multicolumn{1}{l}{\sc Red} & 191 (75.79 \%) & 45 (17.86 \%) & 16 (6.35 \%) & - & - \\ \midrule
		\multicolumn{1}{l}{}  & 1755 & 1366 & 16 & 307 & 734 \\
		%\bottomrule
	\end{tabular}
\end{table}
%\begin{table}[htbp]
%	\small\sf\centering
%	\caption{Proportion of patients assigned to the ED units for each color tag.}\label{tab:triagetag_units_policlinico}
%	\begin{tabular}{lccccc|c}
%		%\toprule
%		%\cline{2-5}
%		%& \multicolumn{2}{|c||}{} & \multicolumn{2}{|c|}{} \\
%		& MU & SU & RA & MIU & OU \\
%		%& \multicolumn{2}{|c||}{} & \multicolumn{2}{|c|}{} \\
%		\midrule			
%		\multicolumn{1}{l}{\sc White} & - & - & -  & 73.44 \%  & 26.56 \% & 100\%\\ %\hline
%		\multicolumn{1}{l}{\sc Green} & 13.75 \% & 34.81\% & - & 14.41 \% & 37.03\% & 100\%  \\ %\hline
%		\multicolumn{1}{l}{\sc Yellow} & 63.95 \% & 33.67\% & - & - & 2.38 \% & 100\%  \\ %\hline
%		\multicolumn{1}{l}{\sc Red} & 75.79 \% & 17.86\% & 6.35 \% & - & - & 100\% \\ 
%		%\bottomrule
%	\end{tabular}
%\end{table}
\begin{table}[htbp]
	\sf\centering
	\caption{Number (and percentage) of green-tagged patients assigned to MU, SU, and MIU in the daytime (8.00 a.m. -- 8.00 p.m.) and at night (8.00 p.m. -- 8.00 a.m.).}\label{tab:triagetag_units_count_policlinico_V}
	\begin{tabular}{lccc}
		%\toprule
		%\cline{2-5}
		%& \multicolumn{2}{|c||}{} & \multicolumn{2}{|c|}{} \\
		& MU & SU & MIU \\
		%& \multicolumn{2}{|c||}{} & \multicolumn{2}{|c|}{} \\
		\midrule
		\multicolumn{1}{l}{\sc Daytime} & 132 (16.60 \%) & 403 (50.69 \%) & 260 (32.70 \%) \\ %\hline
		\multicolumn{1}{l}{\sc Night} & 116 (34.02 \%) & 225 (65.98 \%) & -  \\  \midrule
		\multicolumn{1}{l}{}  & 248 & 628 & 260 \\
		%\bottomrule
	\end{tabular}
\end{table}
%\begin{table}[htbp]
%	\small\sf\centering
%	\caption{Proportion of green-tagged patients assigned to MU, SU, and MIU in the daytime (8.00 a.m. -- 8.00 p.m.) and at night (8.00 p.m. -- 8.00 a.m.).}\label{tab:triagetag_units_policlinico_V}
%	\begin{tabular}{lccc|c}
%		%\toprule
%		%\cline{2-5}
%		%& \multicolumn{2}{|c||}{} & \multicolumn{2}{|c|}{} \\
%		& MU & SU & MIU \\
%		%& \multicolumn{2}{|c||}{} & \multicolumn{2}{|c|}{} \\
%		\midrule
%		\multicolumn{1}{l}{\sc Daytime} & 16.60 \% & 50.69\% & 32.70 \% & 100\%  \\ %\hline
%		\multicolumn{1}{l}{\sc Night} & 34.02 \% & 65.98\% & - & 100\% \\ 
%		%\bottomrule
%	\end{tabular}
%\end{table}

The percentage of patients requiring more than one medical visit is reported in Table~\ref{tab:timesvisit_policlinico} for each color tag and ED unit. Since the number of patients needing more than visit is small, in the DES model all the entities are assumed to undergo one medical visit and, accordingly, each patient flow is represented sequentially. The probability distributions used to generate values for the final waiting time before the discharge are reported in Table~\ref{tab:timesvisit_policlinico} as well (note that the base unit for time adopted in the DES model is the hour). 
\begin{table}[t]
	\caption{Percentage of patients requiring more than one visit and probability distribution of the final waiting time before discharge for each color tag and ED unit.}\label{tab:timesvisit_policlinico}
	\sf\centering
	\begin{tabular}{llcc}
%		\toprule
		&& MORE THAN ONE VISIT & FINAL WAITING TIME \\ 
		\midrule
		{\sc White} & MIU &  0\% & Weib(1.12, 0.41)\\ \midrule
		& MIU & 1.61\% & Weib(0.83, 0.57)\\
		{\sc Green} & SU & 2.44\% & Weib(0.61, 1.07)\\
		& MU & 4.8\% & Weib(0.57, 4.91)\\ \midrule
		\multirow{ 2}{*}{{\sc Yellow}} & SU & 2.04\% &  Weib(0.63, 2.33)\\
		& MU & 3.86\% &  Weib(0.62, 7.22)\\ \midrule
		& RA & 6.25\% & Weib(0.67, 9.26)\\
		{\sc Red} & SU & 6.52\% & Weib(0.71, 4.47) \\
		& MU & 7.9\% & Weib(0.69, 6.65) \\
		\bottomrule
	\end{tabular}
\end{table}

\subsection{Model verification/validation and design of experiments}\label{sec:model_verification_validation_DOE}

The aim of the DES model of the ED of Policlinico Umberto I is to provide a reliable tool for assessing the impact of changes to the current settings on the overcrowding level. To this end, verification and validation are important steps to be carefully considered. While standard techniques have been adopted to verify the model, such as debugging and model trace, validation has required a more sophisticated analysis. Indeed, a model calibration approach has been used to achieve an accurate simulation output that well reproduces the real system data.	   
\par
Important KPIs to describe the ED status are the two time differences contained in ${\cal T}$ and the number of entities generated by the simulation and associated with each color tag. While the former KPIs are part of the objective function used in the calibration procedure and they are calculated across the iterations of the optimization algorithm, the latter do not require a continuous monitoring, since their values are not affected by the calibration. Indeed, the only variations are due to the different number of entities discharged before the simulation ends, which leads to different results in terms of KPI computation. Table~\ref{tab:distr_triage_validation} reports the values of the patient counts obtained from the simulation of the starting point of the calibration procedure along with their confidence intervals. By comparing these values with the counts in Table~\ref{tab:triagetag_units_count_policlinico}, it can be observed that the simulation model provides an accurate output in terms of number of patients.
%\begin{table*}[htbp]
%	\small\sf\centering
%	\caption{Output values (with confidence interval) for the number of patients returned by the simulation of the starting point of the calibration procedure.}\label{tab:distr_triage_validation}
%	\centering
%	\begin{tabular}{lcccc}
%		\toprule
%		% after \\: \hline or \cline{col1-col2} \cline{col3-col4} ...
%		& {\sc White} & {\sc Green} & {\sc Yellow} & {\sc Red}\\ \midrule
%		SU & $$ & $609 \pm 8.10$ & $704 \pm 7.72$ & $54 \pm 3.44$ \\
%		MU & $$ & $264 \pm 5.83$ & $1318 \pm 21.51$ & $220 \pm 5.54$\\
%		RA & & & & $20 \pm 1.45$\\
%		MIU & $42\pm2.47$ &  $241\pm6.20$ & & \\
%		\bottomrule
%		%		& $163.09\pm2.73$ & $1573.49\pm7.59$ & $240.10\pm3.35$ & $17.69\pm0.80$\\
%	\end{tabular}
%\end{table*}
\begin{table}[htbp]
	\sf\centering
	\caption{Output values (with confidence interval) for the number of patients returned by the simulation of the starting point of the calibration procedure.}\label{tab:distr_triage_validation}
	\begin{tabular}{lccc}
		%\toprule
		%\cline{2-5}
%		\toprule
		%& \multicolumn{2}{|c||}{} & \multicolumn{2}{|c|}{} \\
		& MU & SU & MIU \\
		%& \multicolumn{2}{|c||}{} & \multicolumn{2}{|c|}{} \\
		\midrule
		\multicolumn{1}{l}{\sc White} & - & -  & $41\pm2.82$ \\ %\hline
		\multicolumn{1}{l}{\sc Green} & $259 \pm 5.33$ & $612 \pm 9.73$ &  $239\pm6.98$  \\ %\hline
		\multicolumn{1}{l}{\sc Yellow} & $1335 \pm 20.55$ & $716 \pm 9.39$ & -   \\ %\hline
		\multicolumn{1}{l}{\sc Red} & $213 \pm 9.39$ & $57 \pm 5.35$ & - \\ 
%		\midrule
%		\bottomrule
	\end{tabular}
\end{table}

Since the objective and constraint functions compare the real system data and the simulation output, achieving adequate accuracy is imperative for a fair comparison. To this end, the simulation output is estimated through 30 independent simulation replications, each of them 38 days long, with a warm up period of 7 days, thus matching the 31 days of January.

\section{Experimental results}\label{sec:model_calibration_policlinico}
The sets defined in Section~\ref{sec:SBO_calibration} can be adapted to the case study as follows.
\begin{itemize}
	\vspace{-0.2cm}
	\item Let $C = \{\text{W,G,Y,R}\}$ be the set of the triage tags.
	\item Let $U(c) \subseteq \{\text{MU},\text{SU},\text{MIU}\}$ be the set of the ED units where patients with tag $c \in C$ can be visited and treated. In particular, 
	\begin{itemize}
		\item $U(W) = \{\text{MIU}\}$,
		\item $U(G) \ = \{\text{MU},\text{SU},\text{MIU}\}$,			
		\item $U(Y) \ = \{\text{MU},\text{SU}\}$,
		\item $U(R) \ = \{\text{MU},\text{SU}\}$.
	\end{itemize}
\end{itemize}
For the specific instance represented by this case study, 8 different pairs of parameters are considered for the Weibull distributions representing the service times of medical visit and exams. Each pair is associated with an element of the sets $U(c)$. Instead, 7 pairs are used for the probability distributions of triage for the reasons stated in Section~\ref{sec:desm_policlinico}, thus resulting in 23 overall pairs of parameters, i.e., 46 unknown parameters to be determined through the calibration procedure. 
%Further parameters treated as decision variables are the upper boundary parameters of the uniform probability distributions introduced to match the setup and idle times that affect the real system data, thus increasing the number of variables to 60 (a total of 8 uniform distributions are considered among the patient flows associated with the four color tags). 
\par
A crucial and preliminary step is the choice of the starting point of the optimization, given by $x^0$, $y^0$, and $z^0$ (where $x$, $y$, and $z$ are the vectors containing all the corresponding shape and scale parameters, as defined in Section~\ref{sec:SBO_calibration}), which appears not to be straightforward since the missing data prevents the computation of good initial values for the parameters. Instead of starting from randomly generated values, a better strategy is to leverage the known information in a reasonable manner. While the parameters $x_{c \, u}$ of the triage probability distributions are arbitrarily fixed so that the generated service times are in accordance with the ED staff suggestions, for the parameters of the medical visit and exams distributions a different approach is followed. Indeed, in addition to the timestamps shown in Figure~\ref{fig:missing_data}, which are the most commonly known timestamps in every ED, for this specific case study further available information is represented by the time at which the physician requires exams. In general, this happens during the medical visit, although sometimes exams are required throughout the patient flow when periodic check-ups are performed. To provide the parameters $y_{c \, u}$ of the medical visit probability distributions with initial values, a good starting point can be obtained by computing for each patient the duration between the start of the visit and the latest time at which an exam is required. The latter time is considered within a reasonable period from the start of the visit that should reflect the medical visit's maximum service time. Once this timestamp is identified and the durations are consequently available, their values can be used to initialize the parameters by fitting the corresponding Weibull probability distributions. Moreover, the time between the presumptive end of the visit and the discharge can be used in the same manner to obtain good initial values for the parameters $z_{c \, u}$ of the exams distributions.

Although Problem~\ref{prob:model_calibration} is a continuous optimization problem, its variables are considered as discrete to efficiently solve the problem. In particular, at the time of solving the problem, $x^{c \, u}_p$, $y^{c \, u}_p$, and $z^{c \, u}_p$, with $p \in \{1,2\}$, are treated as granular variables, i.e., variables with a controlled number of decimals (see, e.g., \cite{audet.2019}). This choice is motivated by the fact that the output of the simulation model is insensitive to small changes in the values of the decision variables, thus making discrete variables preferable. To this end, a specific granularity is considered based on the role of the parameter associated with the variable in the corresponding probability distribution. All the variables are pairs of shape and scale parameters, which have a different impact on the simulation output. Indeed, while the former parameters determine the shape of the probability distribution, the latter affect the scale of the values generated. With respect to the Weibull distribution adopted, in practice the behavior observed is that the larger the value of the shape parameter, the larger the dispersion of the values of the random variables associated. Moreover, such random variables take on values with a larger scale as the value of the scale parameter increases. Since some preliminary analyses have shown that the impact on the simulation output, measured in terms of the two KPIs $DOT$ and $DIT$, is more noticeable for variations in the scale parameter, for each pair of variables the granularity $\delta^{min}_p$ of the values is fixed to $10^{-3}$ if $p=1$ (i.e., shape parameter) and to $10^{-4}$ if $p=2$ (i.e., scale parameter). Note that treating the variables as granular requires the new variables to be equal to $x^{c \, u}_p/\delta^{min}_p\in \mathbb{Z}$, thus meaning that $x^{c \, u}_p$ have to take on real values that are multiple of $\delta^{min}_p$. The same type of constraint applies to the variables of the visit and exams probability distributions.
\par
To solve the SBO problem for calibrating the simulation model of the ED in hand, the approach of the Sample Average Approximation (SAA) \cite{pagnoncelli.2009,kim.2015-b} is adopted. As a consequence, the resulting optimization problem is deterministic and it can be solved by applying an algorithm from the class of DFO \cite{conn.2009,audet.2017}. In particular, the optimization algorithm proposed in \cite{liuzzi.2020} is used for solving the problem by adopting its default parameters. This algorithm has been successfully applied to the same case study in \cite{desantis.2021} and it is suited for efficiently solving integer black-box constrained problems, like Problem~\ref{prob:model_calibration}, thanks to unconventional search directions, an effective penalty approach, and strong global convergence properties. The maximum number of function evaluations, which represents the stopping condition, is set to 3000. Note that by using the SAA approach, the empirical cumulative distribution functions used in the optimization problem are estimated through the corresponding sample means over the 30 independent simulation replications.
\par
As concerns the optimization problem to solve, the values of the lower and upper bounds introduced in Problem~\eqref{prob:model_calibration} are reported in Table~\ref{tab:lower_upper_bound_calibration}. Moreover, the tolerance $\tau^{c \, u \, i}_{\mu}$ and $\tau^{c \, u \, i}_{\sigma}$ are both fixed to $0.35$ for yellow and green-tagged patients visited and treated in MU and SU, $0.2$ otherwise. 
\begin{table}[htbp]
	\caption{Values of the lower and upper bounds for each pair of variables $(x_1^{c \, u},x_2^{c \, u})$, $(y_1^{c \, u},y_2^{c \, u})$, and $(z_1^{c \, u},z_2^{c \, u})$, for all $c \in C$ and $u \in U(c)$. }\label{tab:lower_upper_bound_calibration}
	\centering		
	\begin{tabular}{cccccc}
		\hline
		\textbf{$l_{x^{c \, u}_1}$} & \textbf{$l_{x^{c \, u}_2}$} & \textbf{$l_{y^{c \, u}_1}$} & \textbf{$l_{y^{c \, u}_2}$} & \textbf{$l_{z^{c \, u}_1}$} & \textbf{$l_{z^{c \, u}_2}$}  \\ \hline
		0.01 & 0.01 & 0.01 & 0.01 & 0.01 & 0.01  \\ \hline\hline
		\textbf{$u_{x^{c \, u}_1}$} & \textbf{$u_{x^{c \, u}_2}$} & \textbf{$u_{y^{c \, u}_1}$} & \textbf{$u_{y^{c \, u}_2}$} & \textbf{$u_{z^{c \, u}_1}$} & \textbf{$u_{z^{c \, u}_2}$}  \\ \hline
		1000 & 0.5 & 1000 & 4 & 1000 & 40  \\ 
		\hline
	\end{tabular}
\end{table} 

Given these settings, the experimental results are obtained by using a PC with Intel Core i7-4790K quad-core 4.00 GHz Processor and 32 GB RAM. Table~\ref{tab:results} reports the service time distributions of triage, medical visit, and exams and reassessments at the final solution determined by the optimization algorithm (we recall that the base unit for time adopted in the DES model is the hour).
\begin{table}[t]
	\caption{Service time distributions of triage, medical visit, and exams and reassessments at the final solution determined by the optimization algorithm. Note that the triage service time of green-tagged patients in MIU is missing since in the DES model these patients are subject to the Weibull distributions used for SU and MU.}\label{tab:results}
	\sf\centering
	\begin{tabular}{llccc}
		%		\toprule
		&& TRIAGE & MEDICAL VISIT & EXAMS \\ 
		\midrule
		{\sc White} & MIU &  Weib(0.61, 0.42) & Weib(0.99, 0.57) & Weib(0.62, 0.23)\\ \midrule
		& MIU & - & Weib(1.06, 0.41) & Weib(1.29, 1.61)\\
		{\sc Green} & SU & Weib(1000, 0.50) & Weib(0.64, 0.31) & Weib(0.68, 2.11)\\
		& MU & Weib(0.60, 0.49) & Weib(1000, 1.63) & Weib(0.69, 12.26)\\ \midrule
		\multirow{ 2}{*}{{\sc Yellow}} & SU & Weib(1000, 0.28) &  Weib(0.51, 0.22) & Weib(0.72, 6.32)\\
		& MU & Weib(0.55, 0.50) &  Weib(0.62, 7.22) & Weib(0.78, 25.46)\\ \midrule
		\multirow{ 2}{*}{{\sc Red}} & SU & Weib(0.92, 0.10) & Weib(0.88, 0.75) & Weib(1.41, 19.34) \\
		& MU & Weib(0.80, 0.19) & Weib(2.24, 0.36) & Weib(0.78, 30.92) \\
		\bottomrule
	\end{tabular}
\end{table}
Two types of plots are shown in Appendices \ref{app:appB} and \ref{app:appC} to assess the accuracy of the calibrated simulation model.
% Figures~\ref{fig:isto_red_SU_TA}--\ref{fig:isto_white_MIU_TB} 
The former appendix focuses on the KPI given by the average hourly number of patients inside the ED after the start of the medical visit, namely, patients satisfying two conditions: $(i)$ they have already started the medical visit; $(ii)$ they have not been discharged yet. In particular, the plots report the comparison between the value of the KPI computed through the real data and the corresponding value derived from the simulation output along with its $95\%$ confidence interval.
%Moreover, Figures~\ref{fig:ecdf_red_SU_TA}--\ref{fig:ecdf_white_MIU_TB} 
The latter appendix reports the comparison between the Empirical Cumulative Distribution Functions (ECDFs) of the values of the time differences $DOT$ and $DIT$ computed for each patient through the real data and the simulation output resulting from the model calibration procedure. In both appendices, the comparisons are performed for all the color tags $c \in C$ and for all the ED units $u \in U(c)$ where a patient with tag $c$ can be assigned. Moreover, both the average hourly number of patients inside the ED after the start of the medical visit and the ECDFs are obtained as an average over the independent simulation replications. Both types of plots show that, in general, the values of the KPIs corresponding to the system data are well reproduced by the simulation output. However, some dissimilarities may be still present in spite of the calibration, especially for green and yellow-tagged patients. This is due to the difficulty in reproducing the patient flow when sharing of the resources is involved between patients with different triage tag, as is the case for the medical visit of green and yellow-tagged patients at both $MU$ and $SU$. Moreover, the experimental results show that, in general, the comparison between real data and simulation output is more accurate for the time differences DIT. Indeed, larger errors are observed in both types of plots when time differences DOT are considered due to the difficulty in reproducing the actual waiting times of patients, which depend both on the service time of the medical visit and on the number of patients that are in queue. Hence, it is the result of a {\em seize-delay-release} process used within the simulation, while DIT is given by the mere sum of service times. 

By observing the ECDF plots in Appendix \ref{app:appC}, it is possible to gain useful insight that histograms, which are commonly adopted in the literature, fail to provide. For instance, the ECDF related to yellow-tagged patients in SU shows that the time differences DOT collected from the simulation model are systematically larger than the corresponding real values, meaning that either the estimated parameters of the probability distribution of the medical visit service time or the modeling of the interaction with green-tagged patients in SU may be improved with further research. To avoid ambiguities, it is important to remark that since the simulated ECDFs are obtained as an average over independent simulation replications, the jumps (or steps) in these simulated functions correspond to patients observed across all the replications. This is a consequence of averaging different step functions, such as the ECDFs, which have jumps at different points in their domain.  
\par
\begin{figure}[htbp]
	\centering
	\includegraphics[width=\textwidth]{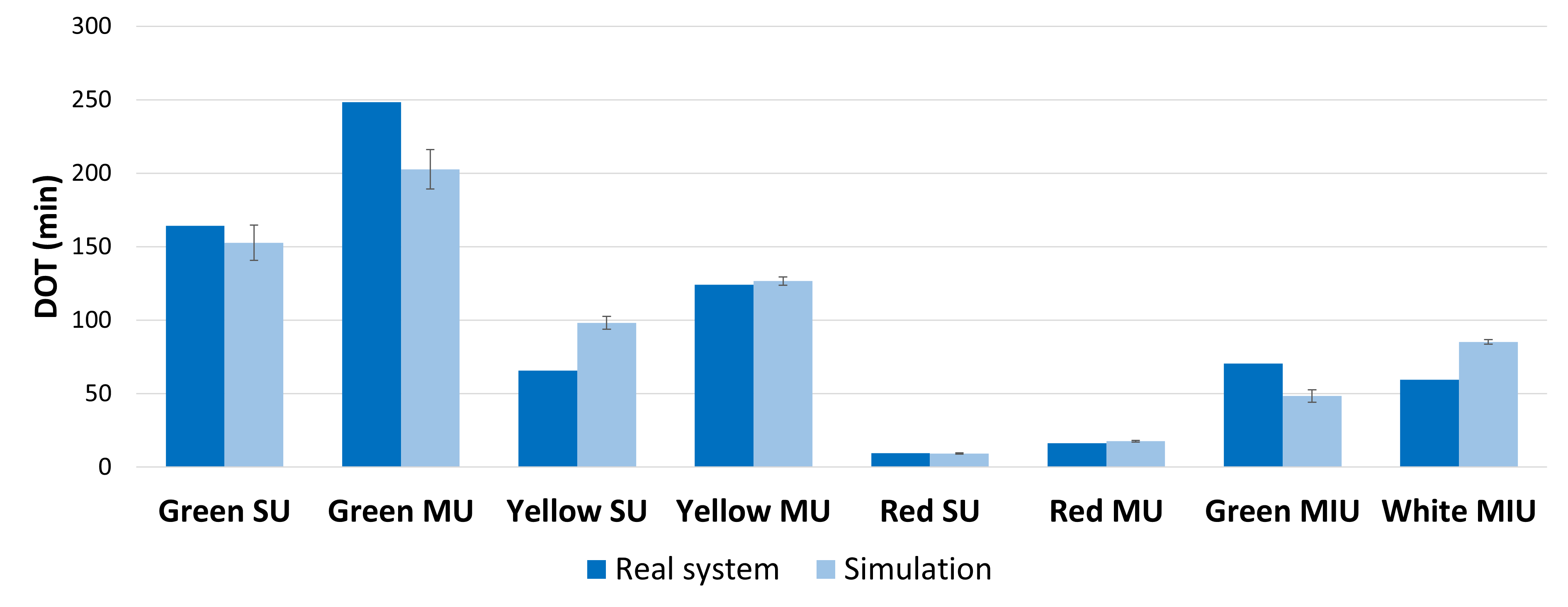}
	\caption{Plot of current values and simulation
		output of DOT with the confidence interval.}\label{fig:confronto_medie_TA_calibration}
\end{figure}
\begin{figure}[htbp]
	\centering
	\includegraphics[width=\textwidth]{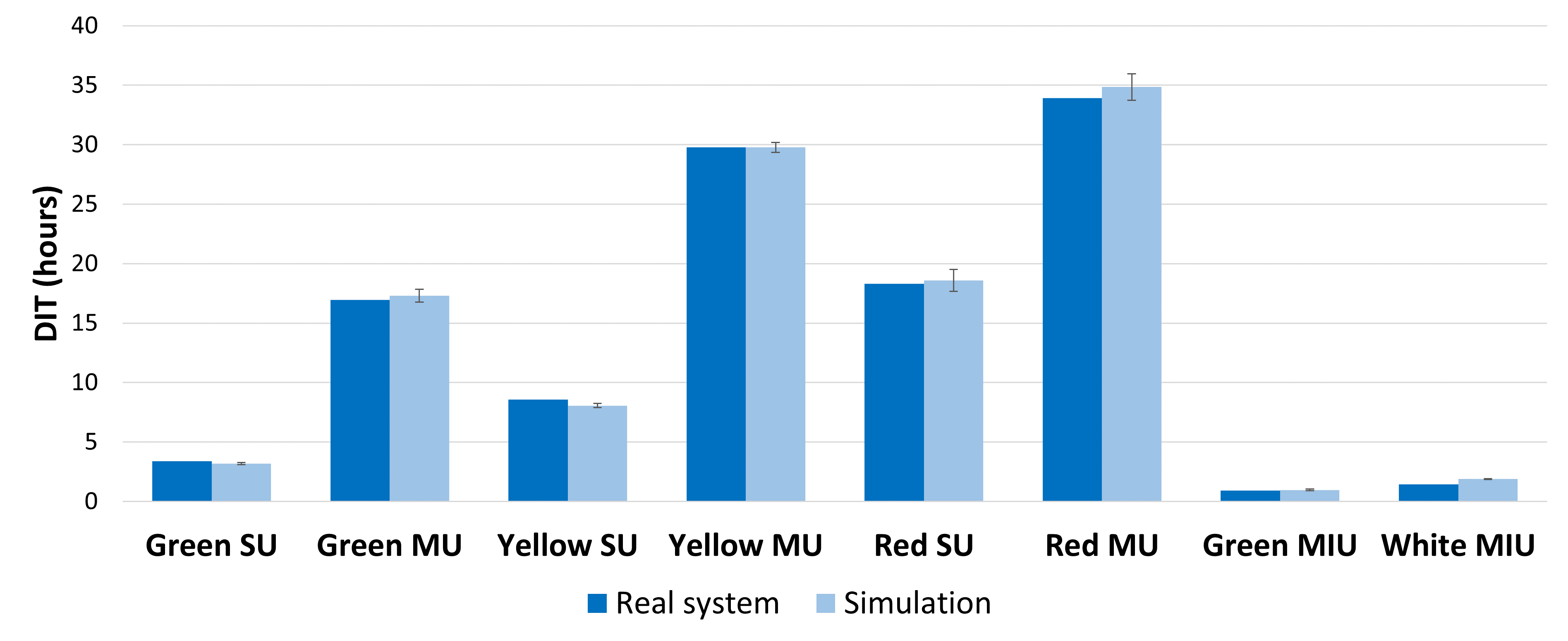}
	\caption{Plot of current values and simulation
		output of DIT with the confidence interval.}\label{fig:confronto_medie_TB_calibration}
\end{figure}

The same conclusions on the accuracy of the simulation model can be drawn also by observing Figures~\ref{fig:confronto_medie_TA_calibration}--\ref{fig:confronto_medie_TB_calibration}, which report the comparison between the average current and simulated values of DOT and DIT for all the triage tags and all the ED units.  As already mentioned, the simulation output related to the time differences DIT is a good approximation of the real
system values since for each color tag the current values are
either within the corresponding confidence interval or close
to it. Instead, even though for the time differences DOT the accuracy of the simulation responses as estimates of the current values is lower, the relative error is within the tolerances $\tau^{c \, u \, i}_{\mu}$ chosen. Indeed, all the constraints of Problem \eqref{prob:model_calibration} are (of course) satisfied at the optimal solution obtained by applying the model calibration approach. Note that reducing the values of the tolerances potentially allows achieving more accurate solutions, however, the optimization algorithm is not guaranteed to determine a feasible solution.    
\par
The experimental results discussed above show that the model calibration procedure enables estimating the missing parameters in order for the simulation model to satisfactorily reproduce the ED operations despite the unavailable timestamps, which represent the main hurdle for achieving a high level of accuracy. The constraints on the sample means used in Problem \eqref{prob:model_calibration} ensure that the responses of the simulation model are on average close to the corresponding real values. Moreover, the constraints on the sample variances of the KPIs help the optimization algorithm to avoid points associated with inaccurate results. However, when looking more closely into the details of the numerical results, some dissimilarities between real data and simulation output may be observed as a result of the impact of the problem of missing data, which necessarily undermines the overall accuracy. Notwithstanding, the simulation model may be deemed sufficiently reliable with respect to the specific objectives of the analysis. 

\section{Conclusions}\label{sec:conc}
The SBO approach proposed in this paper addresses one of the data quality problems that frequently affect the DES models concerning EDs. In particular, the focus is on the problem of missing data, which consists in the unavailability of data related to some of the starting and ending times of the activities performed in the ED. This well-known issue is responsible for the lack of knowledge of the service time of some processes, which is required to estimate the corresponding probability distributions to use in the simulation model. For this purpose, after assuming suitable probability distributions for the processes associated with the unknown service times, a model calibration procedure is proposed to estimate the missing parameters of such distributions. 

The preliminary experimental results on the ED of a big hospital in Italy show that the model calibration procedure enables estimating the missing parameters in order for the simulation model to satisfactorily reproduce the ED operations despite the unavailable timestamps, which represent the main hurdle for achieving a high level of accuracy. Although proper constraints ensure that the responses of the simulation model are on average close to the corresponding real values, when looking more closely into the details of the experimental results, some dissimilarities between real data and simulation output may be observed as a result of the impact of the problem of missing data, which necessarily undermines the overall accuracy. Notwithstanding, the simulation model may be deemed sufficiently reliable with respect to the specific objectives of the analysis. 
Since there is still a margin for improvement, several ideas are considered for further research. For instance, as the worst results are observed for the time differences DOT, a significant improvement may be obtained by using different weights for the terms of the objective function of Problem \eqref{prob:model_calibration} and by assigning larger values to the terms associated with DOT. Moreover, different objective functions and starting points may be taken into account also to assess the robustness of the approach. Finally, increasing the number of function evaluations used as stopping condition may lead to better solutions allowing a more thorough exploration of the feasible region.

\section*{Acknowledgment}
The authors are grateful to Prof. F. Romano (Chief Medical Officer) and Dr. L. De Vito (Medical Director of ED) of \textit{Policlinico Umberto I} of Rome for their availability to carry out this study.

\appendix

\section{Model calibration - Plots I}
\label{app:appB}
This appendix focuses on the KPI given by the average hourly number of patients satisfying two conditions: $(i)$ they have already started the medical visit; $(ii)$ they have not been discharged yet. In particular, the plots report the comparison between the value of the KPI computed through the real data and the corresponding value derived from the simulation output along with its $95\%$ confidence interval. The comparisons are performed for all color tags $c \in C$ and for all the ED units $u \in U(c)$ where a patient with color tag $c$ can be assigned.
%\begin{figure}[htbp]		
%	\label{fig:isto_red_SU_TA}
%	\centering
%	\includegraphics[width=0.65\textwidth]{red_SU_TA.png} 	
%\end{figure}
%\begin{figure}[htbp]
\begin{figure}[H]
	\centering
	\scalebox{.8}{\includegraphics[width=\textwidth]{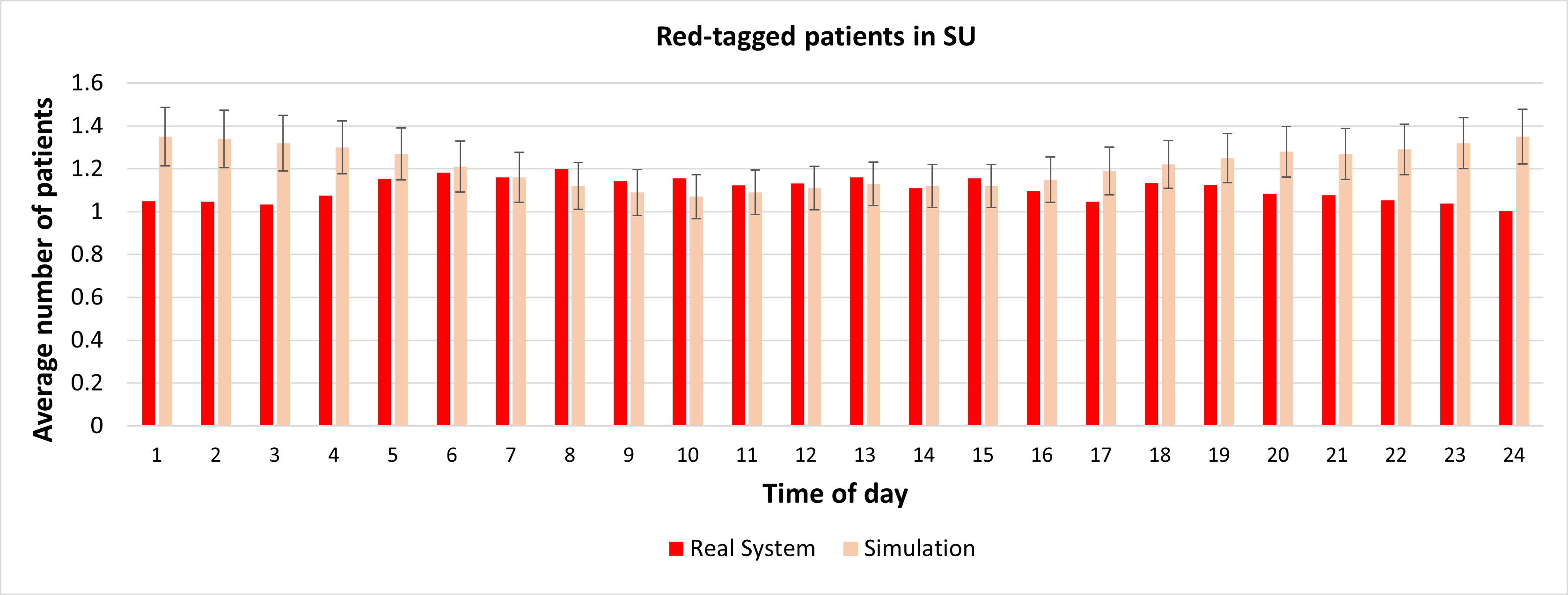}}
	\label{fig:isto_red_SU_TB}
\end{figure}
%\begin{figure}[htbp]
%	\centering
%	\includegraphics[width=0.65\textwidth]{red_MU_TA.png} 
%	\label{fig:isto_red_MU_TA}
%\end{figure}
%\begin{figure}[htbp]
\begin{figure}[H]
	\centering
	\scalebox{.8}{\includegraphics[width=\textwidth]{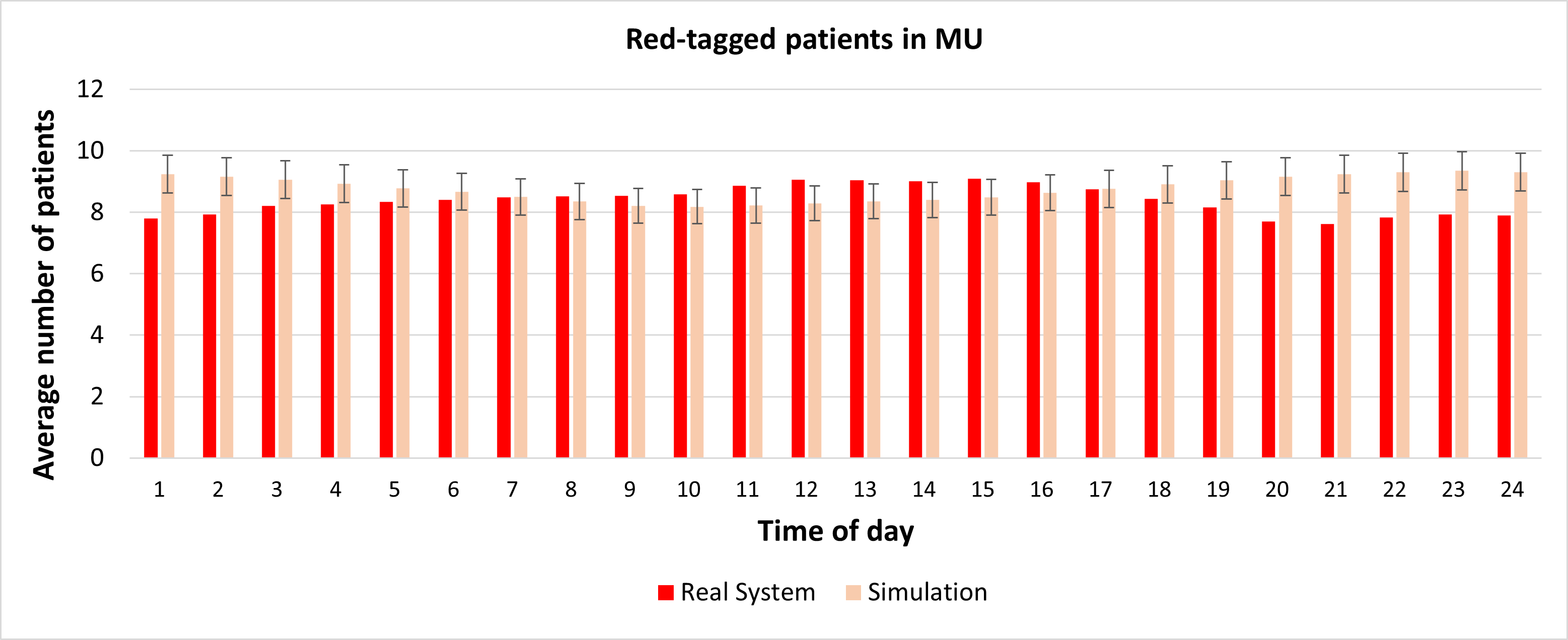}}
	\label{fig:isto_red_MU_TB}
\end{figure}
%\begin{figure}[htbp]
%	\centering
%	\includegraphics[width=0.65\textwidth]{red_RA_TA.png} 
%	\label{fig:isto_red_RA_TA}
%\end{figure}
%\begin{figure}[htbp]
%	\centering
%	\includegraphics[width=\textwidth]{red_RA_TB.png}
%	\label{fig:isto_red_RA_TB}
%\end{figure}
%\begin{figure}[htbp]
%	\centering
%	\includegraphics[width=0.65\textwidth]{yellow_SU_TA.png} 
%	\label{fig:isto_yellow_SU_TA}
%\end{figure}
%\begin{figure}[htbp]
\begin{figure}[H]
	\centering
	\scalebox{.8}{\includegraphics[width=\textwidth]{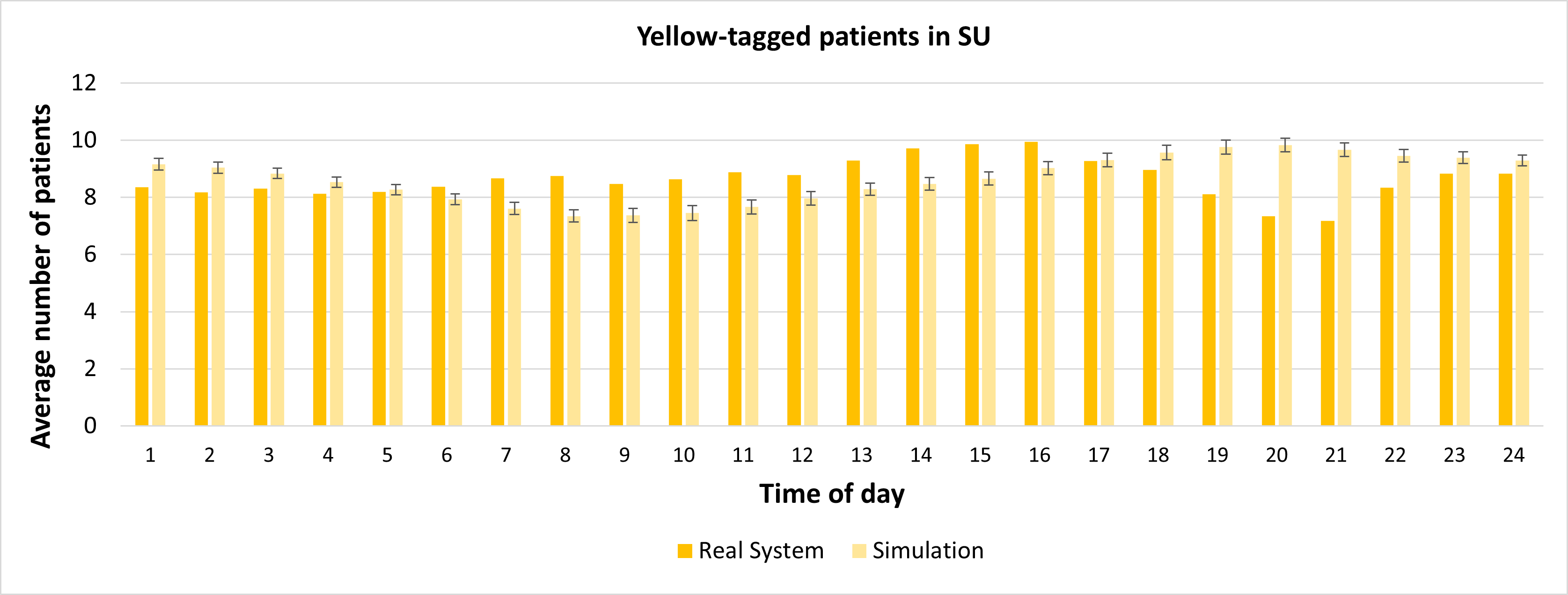}}
	\label{fig:isto_yellow_SU_TB}
\end{figure}
%\begin{figure}[htbp]
%	\centering
%	\includegraphics[width=0.65\textwidth]{yellow_MU_TA.png} 
%	\label{fig:isto_yellow_MU_TA}
%\end{figure}
%\begin{figure}[htbp]
\begin{figure}[H]
	\centering
	\scalebox{.8}{\includegraphics[width=\textwidth]{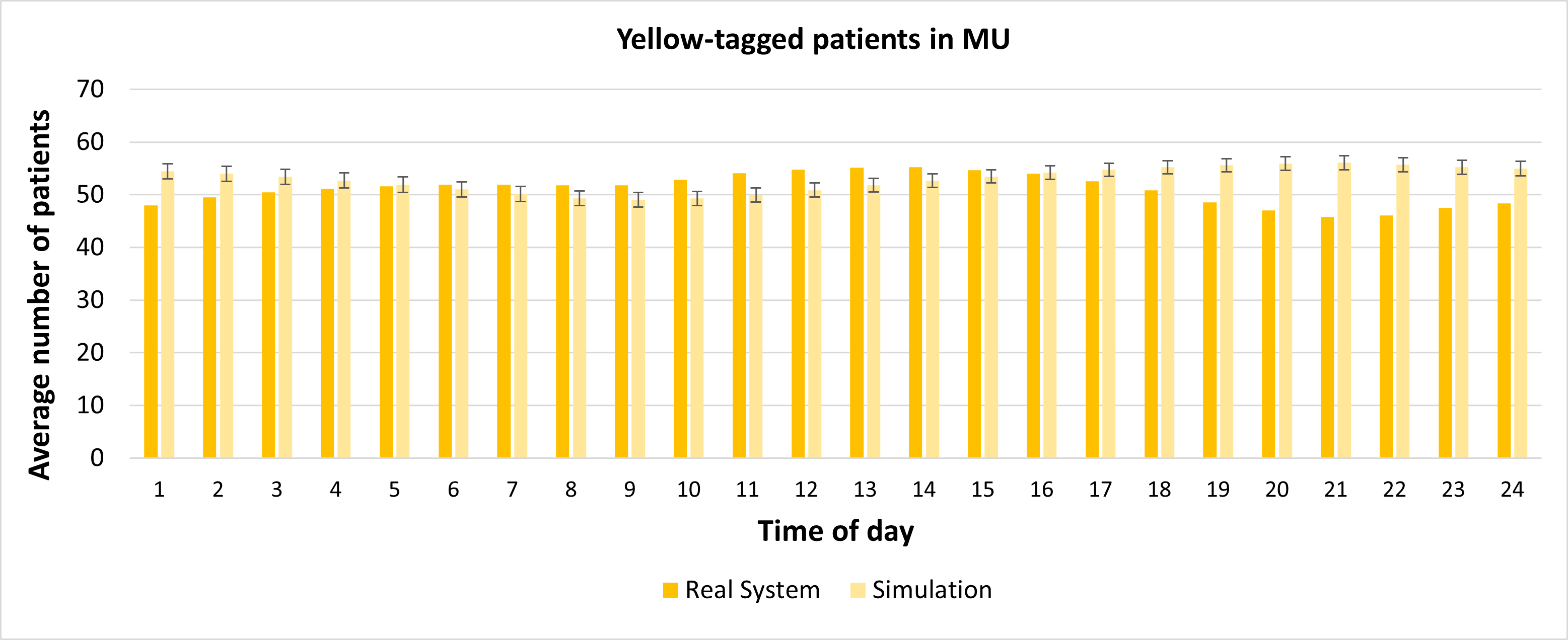}}
	\label{fig:isto_yellow_MU_TB}
\end{figure}
%\begin{figure}[htbp]
%	\centering
%	\includegraphics[width=0.65\textwidth]{green_SU_TA.png} 
%	\label{fig:isto_green_SU_TA}
%\end{figure}
%\begin{figure}[htbp]
\begin{figure}[H]
	\centering
	\scalebox{.8}{\includegraphics[width=\textwidth]{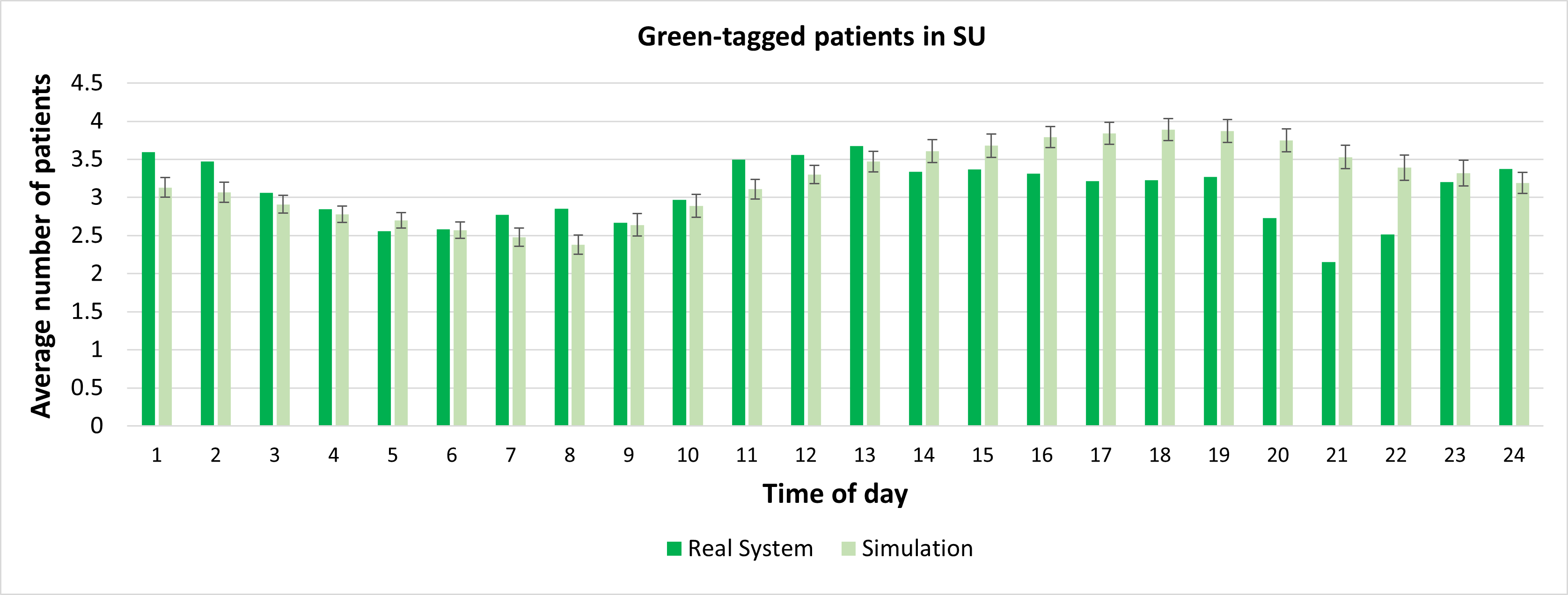}}
	\label{fig:isto_green_SU_TB}
\end{figure}
%\begin{figure}[htbp]
\begin{figure}[H]
%	\centering
%	\includegraphics[width=0.65\textwidth]{green_MU_TA.png} 
%	\label{fig:isto_green_MU_TA}
%\end{figure}
%\begin{figure}[htbp]
	\centering
	\scalebox{.8}{\includegraphics[width=\textwidth]{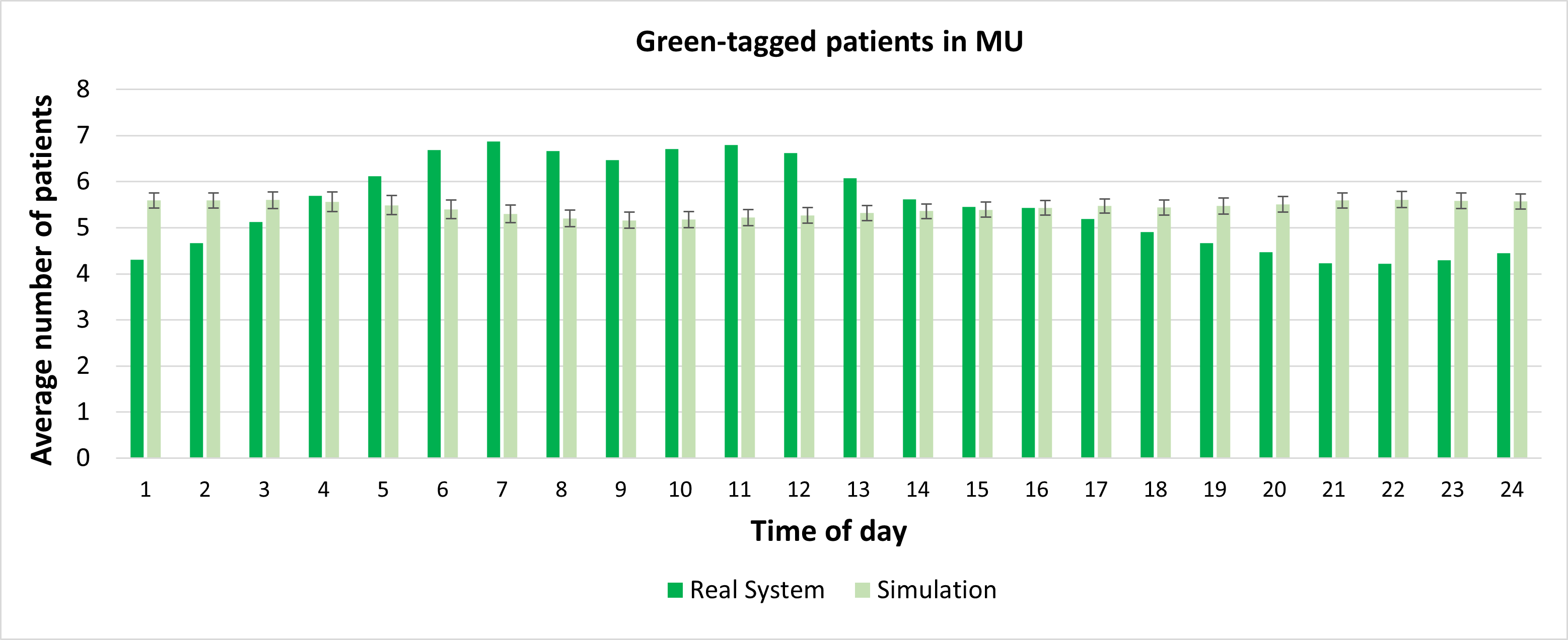}}
	\label{fig:isto_green_MU_TB}
\end{figure}
%\begin{figure}[htbp]
%	\centering
%	\includegraphics[width=0.65\textwidth]{green_MIU_TA.png} 
%	\label{fig:isto_green_MIU_TA}
%\end{figure}
%\begin{figure}[htbp]
\begin{figure}[H]
	\centering
	\scalebox{.8}{\includegraphics[width=\textwidth]{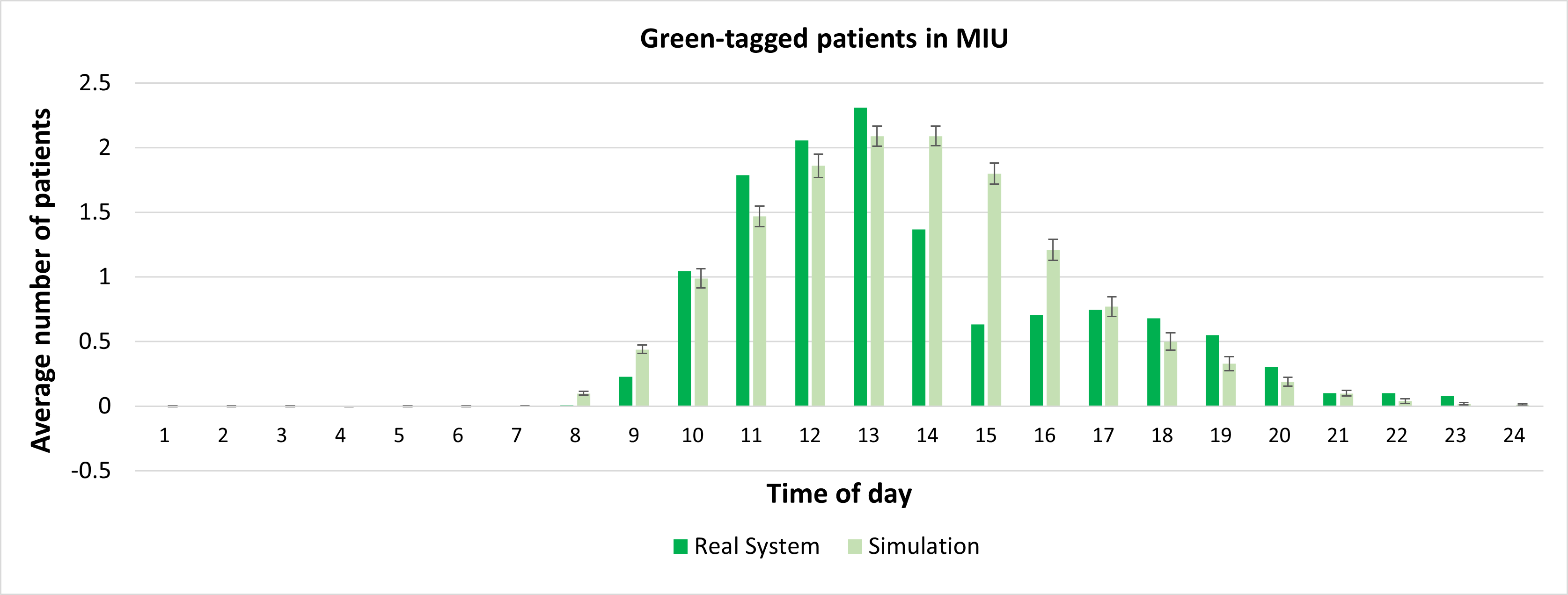}}
	\label{fig:isto_green_MIU_TB}
\end{figure}
%\begin{figure}[htbp]
%	\centering
%	\includegraphics[width=0.65\textwidth]{white_MIU_TA.png} 
%	\label{fig:isto_white_MIU_TA}
%\end{figure}
%\begin{figure}[htbp]
\begin{figure}[H]
	\centering
	\scalebox{.8}{\includegraphics[width=\textwidth]{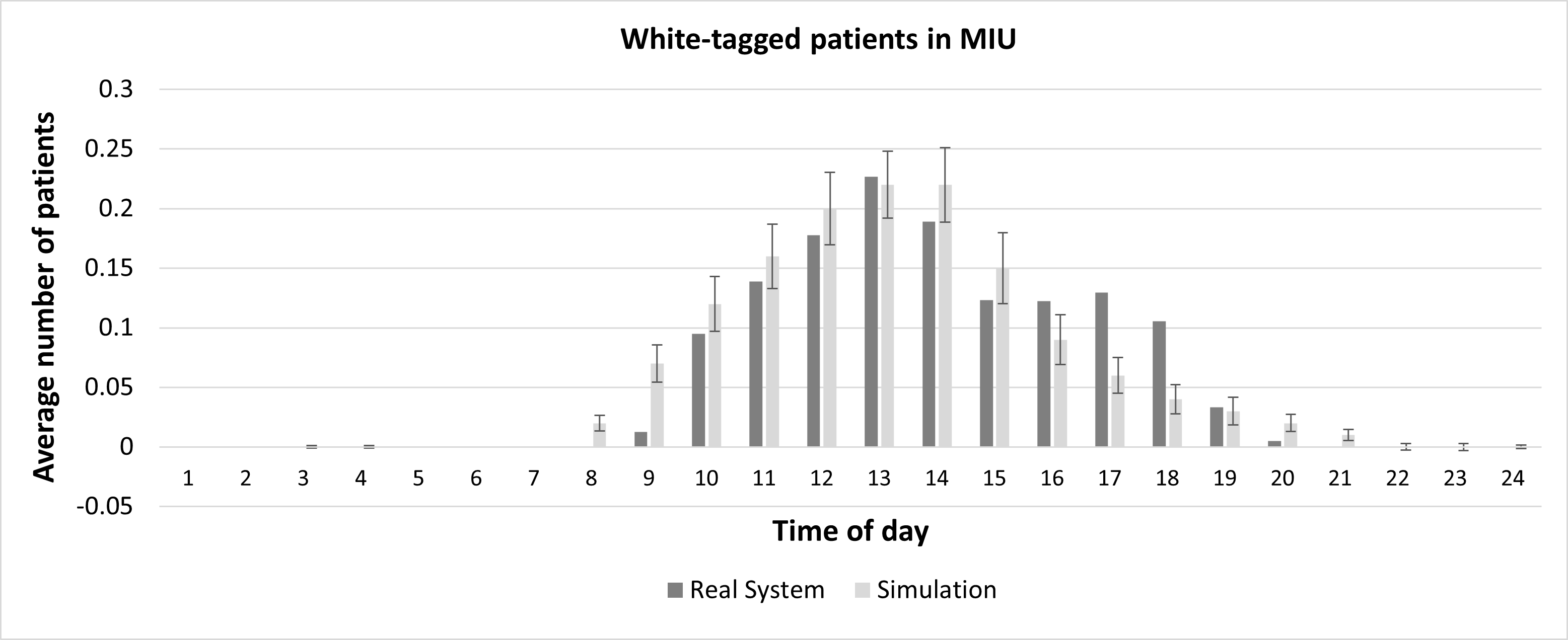}}
	\label{fig:isto_white_MIU_TB}
\end{figure}

%\section{Appendix}
\section{Model calibration - Plots II}
\label{app:appC}
This appendix reports the comparison between the Empirical Cumulative Distribution Functions (ECDFs) of the time differences $DOT$ and $DIT$ computed through the real data and the simulation output resulting from the calibration procedure. The comparisons are performed for all color tags $c \in C$ and for all ED units $u \in U(c)$. The colored curves correspond to the simulation output.
%\begin{figure}[htbp]
\begin{figure}[H]
	\centering
	\scalebox{.8}{\includegraphics[width=0.48\textwidth]{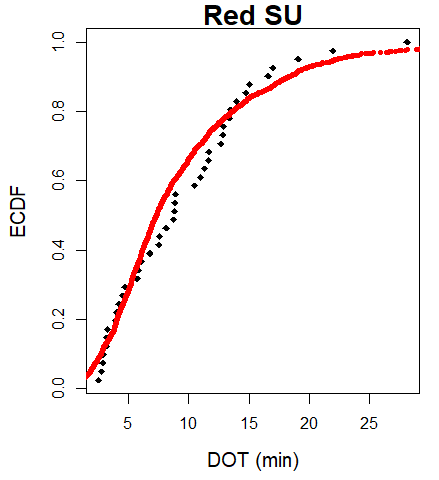}} \quad \scalebox{.8}{\includegraphics[width=0.48\textwidth]{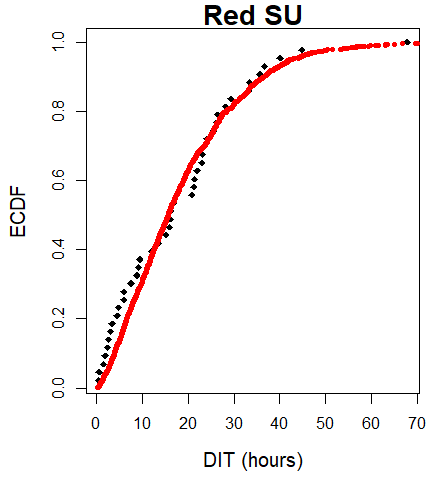}}
	\label{fig:ecdf_red_SU_TAB}
\end{figure}
%\begin{figure}[htbp]
\begin{figure}[H]
	\centering
	\scalebox{.8}{\includegraphics[width=0.48\textwidth]{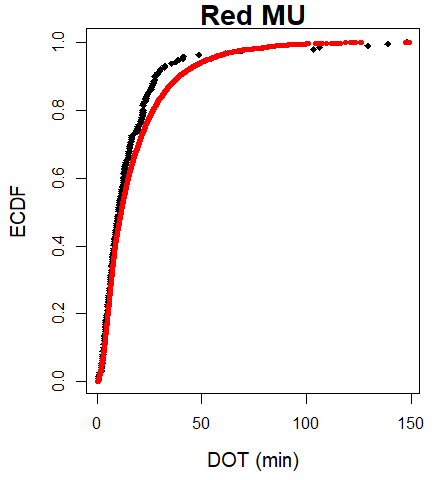}} \quad \scalebox{.8}{\includegraphics[width=0.48\textwidth]{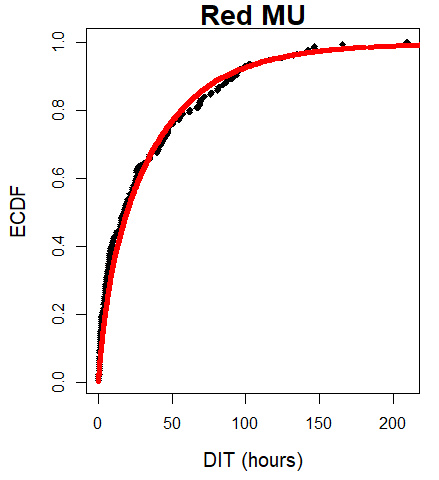}}
	\label{fig:ecdf_red_MU_TAB}
\end{figure}
%\begin{figure}[htbp]
	\begin{figure}[H]
	\centering
	\scalebox{.8}{\includegraphics[width=0.48\textwidth]{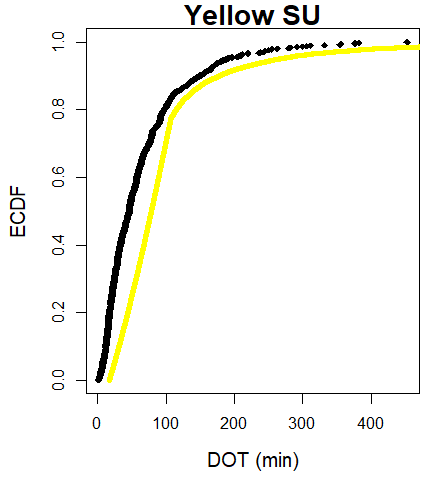}} \quad \scalebox{.8}{\includegraphics[width=0.48\textwidth]{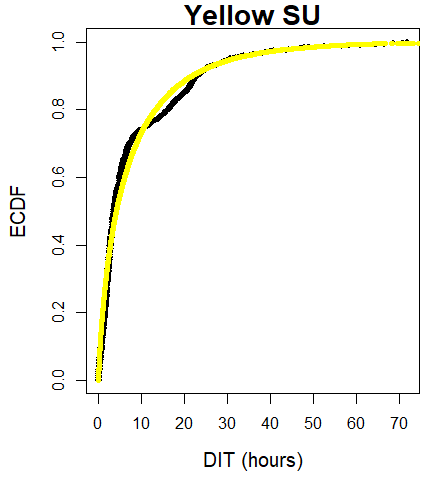}}
	\label{fig:ecdf_red_RA_TAB}
\end{figure}
%\begin{figure}[htbp]
	\begin{figure}[H]
	\centering
	\scalebox{.8}{\includegraphics[width=0.48\textwidth]{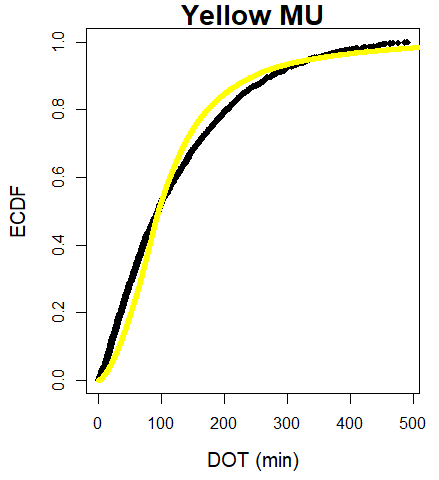}} \quad \scalebox{.8}{\includegraphics[width=0.48\textwidth]{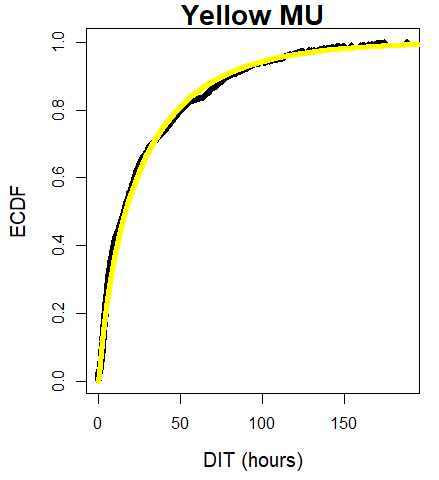}}
	\label{fig:ecdf_yellow_SU_TAB}
\end{figure}
%\begin{figure}[htbp]
	\begin{figure}[H]
	\centering
	\scalebox{.8}{\includegraphics[width=0.48\textwidth]{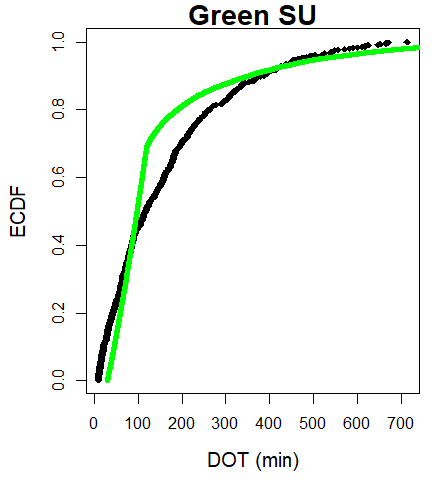}} \quad \scalebox{.8}{\includegraphics[width=0.48\textwidth]{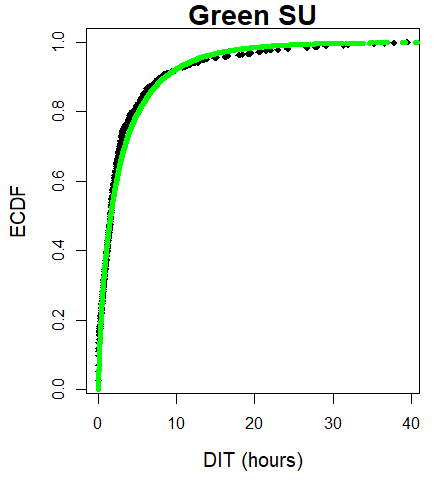}}
	\label{fig:ecdf_yellow_MU_TAB}
\end{figure}
%\begin{figure}[htbp]
	\begin{figure}[H]
	\centering
	\scalebox{.8}{\includegraphics[width=0.48\textwidth]{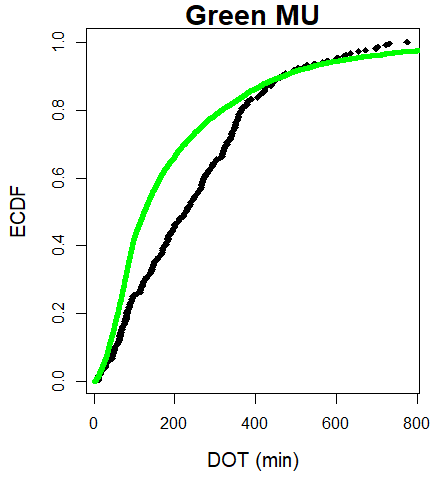}} \quad \scalebox{.8}{\includegraphics[width=0.48\textwidth]{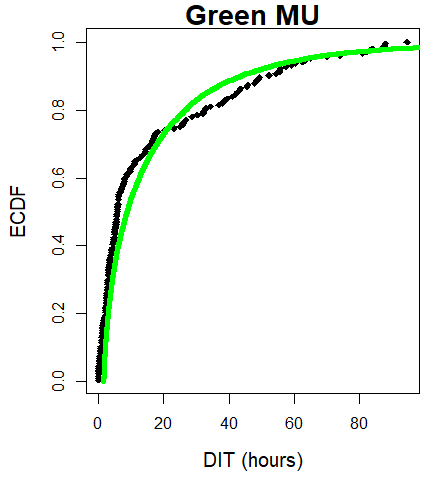}}
	\label{fig:ecdf_green_SU_TAB}
\end{figure}
%\begin{figure}[htbp]
	\begin{figure}[H]
	\centering
	\scalebox{.8}{\includegraphics[width=0.48\textwidth]{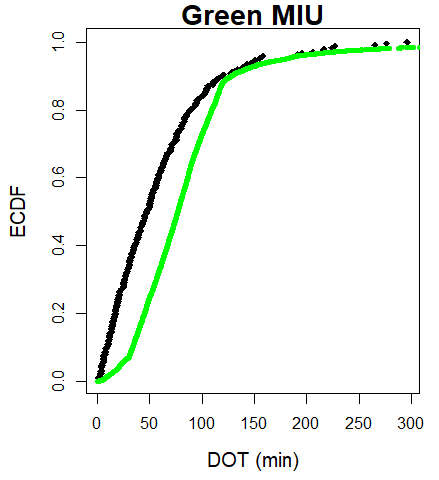}} \quad \scalebox{.8}{\includegraphics[width=0.48\textwidth]{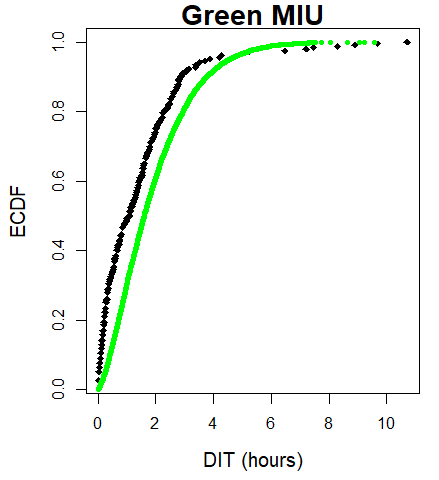}}
	\label{fig:ecdf_green_MU_TAB}
\end{figure}
%\begin{figure}[htbp]
	\begin{figure}[H]
	\centering
	\scalebox{.8}{\includegraphics[width=0.48\textwidth]{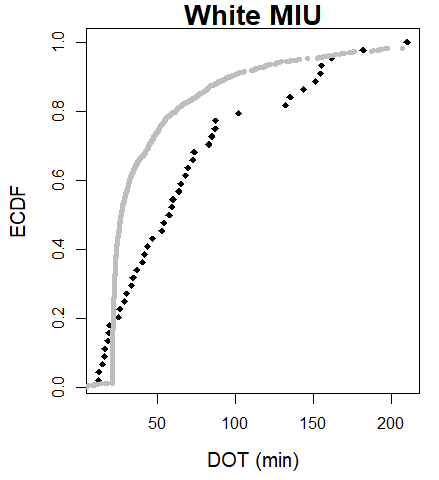}} \quad \scalebox{.8}{\includegraphics[width=0.48\textwidth]{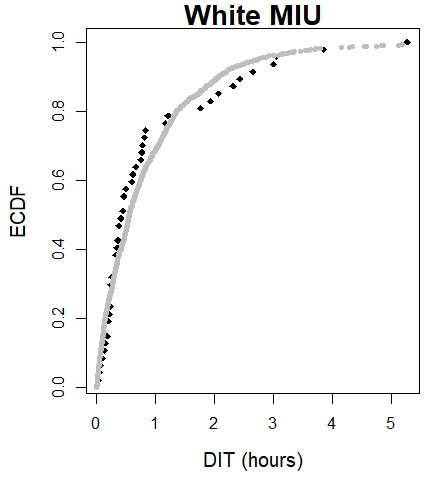}}
	\label{fig:ecdf_green_MIU_TAB}
\end{figure}

%\begin{acknowledgements}
%\end{acknowledgements}

\bibliographystyle{spbasic}      % basic style, author-year citations

%\begin{thebibliography}{}
%%
%% and use \bibitem to create references. Consult the Instructions
%% for authors for reference list style.
%\end{thebibliography}

\section*{Conflict of interest}

The authors declare that they have no conflict of interest.

%\bibliographystyle{elsarticle-num}
%\bibliography{EDreference,simulation,derivativeFree,New_healthcare_reference}

\end{document}